%% file: arxiv0429.tex
\documentclass[12pt]{article}
\UseRawInputEncoding

{\begin{Sbox}\begin{minipage}}%
{\end{minipage}\end{Sbox}\fbox{\TheSbox}}%

\usepackage[psamsfonts]{amssymb}
\usepackage{amsmath, amsthm, anysize, enumerate, color, url}
\usepackage{graphicx}
\usepackage{epstopdf}
\usepackage{epsfig}
\usepackage{subfigure}

\usepackage{natbib}
\usepackage[colorlinks, breaklinks,
                   linkcolor=red,
                  citecolor=blue
                  ]{hyperref}

\usepackage{breakurl}
\usepackage{threeparttable}

\newtheorem{theorem}{Theorem} 
\newtheorem{corollary}{Corollary}

\newtheorem{lemma}{Lemma} 

\theoremstyle{definition}

\newcommand{\hh}{\mathbf{h}}
\newcommand{\dd}{\mathbf{d}}
\newcommand{\E}{\mathbb{E}}

\newcommand{\R}{\mathbb{R}}

\newcommand{\diag}{\mbox{diag}}
\newcommand{\bs}{\boldsymbol}

\usepackage{setspace}

\begin{document}

\title{Wilks's theorems in some exponential random graph models}

\author{Ting Yan\thanks{Department of Statistics, Central China Normal University, Wuhan, 430079, China.
\texttt{Email:} tingyanty@mail.ccnu.edu.cn.}
\hspace{4mm}
Yuanzhang Li\thanks{Walter Reed Army Institute of Research,
503 Robert Grant Ave., Silver Spring, Maryland, 20910, USA.
\texttt{Email:} Liy.Li@us.army.mil}
\hspace{4mm}
Jinfeng Xu\thanks{Department of Statistics and Actuarial Science,
    The University of Hong Kong, Hong Kong, 10016.
\texttt{Email:} xujf@hku.hk}
\hspace{4mm}
Yaning Yang\thanks{Department of Statistics and Finance,
University of Science and Technology of China, Anhui, 230026, China.
\texttt{Email:} ynyang@ustc.edu.cn}
\hspace{4mm}
Ji Zhu\thanks{ Department of Statistics,
University of Michigan, Ann Arbor, Michigan, USA.
\texttt{Email:} jizhu@umich.edu}
}
\date{}

\maketitle

\begin{abstract}
\begin{spacing}{1.2}

We are concerned here with the likelihood ratio statistics in two exponential random graph models--the $\beta$-model and the Bradley--Terry model,
in which the degree sequence on an undirected graph and the out-degree sequence
on a weighted directed graph are the exclusively sufficient statistics in the exponential-family distributions on graphs, respectively.
We prove the Wilks type of theorems for some fixed and growing dimensional hypothesis testing problems.
More specifically,  under two fixed dimensional null hypotheses $H_0: \beta_i=\beta_i^0$ for $i=1,\ldots, r$ and $H_0: \beta_1=\ldots=\beta_r$,
we show that $2[\ell(\widehat{\bs{\beta}}) - \ell(\widehat{\bs{\beta}}^0)]$ converges in distribution to a Chi-square distribution with
the respective degrees of freedoms, $r$ and $r-1$, as the dimension $n$ of the full parameter space goes to infinity. Here,
$\ell( \bs{\beta})$ is the log-likelihood function on the parameter $\bs{\beta}$, $\widehat{\bs{\beta}}$ is the MLE under the full parameter space,
and $\widehat{\bs{\beta}}^0$ is the restricted MLE under the null parameter space.
For two increasing dimensional null hypotheses $H_0: \beta_i = \beta_i^0$ for $i=1, \ldots, n$ and $H_0: \beta_1=\ldots=\beta_r$ with $r/n \ge c$,
we show that the normalized log-likelihood ratio statistics, $(2[\ell(\widehat{\bs{\beta}}) - \ell(\bs{\beta}^0)] -n)/(2n)^{1/2}$ and
$(2[\ell(\widehat{\bs{\beta}}) - \ell(\widehat{\bs{\beta}}^0)] -r)/(2r)^{1/2}$, both
converge in distribution to the standard normal distribution.
Simulation studies and an application to NBA data illustrate the theoretical results.
\end{spacing}

\vskip 5 pt \noindent
\begin{spacing}{1.4}
\end{spacing}

\end{abstract}

\vskip 5pt

\section{Introduction}

We are concerned here with the likelihood ratio tests in two exponential random graph models including
the $\beta$-model [\cite{Chatterjee:Diaconis:Sly:2011}] for undirected graphs and the Bradley--Terry model [\cite{bradley-terry1952}] for weighted directed graphs. These two models are closely related, where each node is assigned one parameter $\beta_i$.
The $\beta$-model postulates that node $i$ is connected to node $j$ with probability $\mu(\beta_i + \beta_j)$
while the Bradley--Terry model assumes that $i$ is preferred to (or wins) $j$ with probability $\mu(\beta_i - \beta_j)$,
where $\mu(x)=e^x/(1+e^x)$ is the logistic function. We will call $\beta_i$ the strength parameter hereafter.

The Bradley--Terry model originated from paired comparison data that can be represented in a weighted directed graph,
where each node denotes one subject and a weighted directed edge  from node $i$ to node $j$ indicating the number of times that subject $i$ is preferred to subject $j$.
``Subject" is a covering term that could stand for team in sports games, journals in citation networks, brand in products and many more.
The Bradley--Terry model has numerous applications ranging from the rankings of 
classical sports teams [\cite{masarotto2012the, Sire_2008baseball, Whelan-2020-Hockey}] and scientific journals [\cite{stigler1994citation, Varin-2016-jrsa}], to the quality of product brands [\cite{radlinski2007active}],
to the transmission/disequilibriumtest in genetics [\cite{Sham:Curtis:1995}] and crowdsourcing [\cite{chen2016overcoming}].
The related $\beta$-model has been widely used to model the degree heterogeneity of random graphs [\cite{Park:Newman:2004,Blitzstein:Diaconis:2011,Chatterjee:Diaconis:Sly:2011}].

Since the number of parameters matches the number of nodes and the sample is only one realized graph,
asymptotic inference is nonstandard and turns to be challenging [\cite{Goldenberg2010,Fienberg2012,Graham2017,Chen:2020}].
There has been received considerable attentions to explore theoretical properties in the $\beta$-model and the Bradley--Terry model
in the past decades. Consistency and asymptotic normality of the maximum likelihood estimator (MLE) are established [\cite{simons-yao1999,Chatterjee:Diaconis:Sly:2011,Yan:Xu:2013,Chen:2020,Han-chen2020}].
We shall elaborate them after stating our results.
One fundamental characteristic of guaranteeing nice asymptotic properties of the MLE
is that $n^2$ observed edges lurks in the graph in contrast with $n$ parameters.
However, little is known about asymptotic properties of the likelihood ratio statistics
for these two models under the high dimensional setting. The aim of this paper is to fill this gap.

The likelihood ratio statistics play a very important role in parametric hypothesis testing problems. 
Under the large sample framework that the dimension of parameters is fixed and the size of samples
goes to infinity, one of the most celebrated results is the Wilks theorem [\cite{wilks1938}].
That says minus twice log-likelihood ratio statistic $-2\log \Lambda(X)$ under the null $H_0:\theta\in\Theta_0$, converges in distribution to
a Chi-square distribution with $k$ degrees of freedom independent of nuisance parameters,
where $ \Lambda(x) =  \max_{\theta\in \Theta_0}  f(\mathbf{x}| \theta) /
 \max_{\theta\in \Theta} f(\mathbf{x}|\theta)$, $f(\mathbf{x}|\theta)$ is a  probability density function for $n$ samples
$\mathbf{x}=(x_1, \ldots, x_n)$ and $k=\mathrm{dim}(\Theta) - \mathrm{dim}(\Theta_0)$.
This appealing property 
was referred to as the Wilks phenomenon by \cite{fan2001}.
Since the dimension of parameter space often increases with the size of samples,
it is interesting to see whether the Wilks type of results continue to hold in the high dimension setting.
In this paper, we investigate the Wilks' theorems for the Bradley--Terry model and the $\beta$-model in both fixed and increasing dimensional null hypotheses.
Our contributions are as follows.
\begin{itemize}
\item  Under two fixed dimensional null hypotheses $H_0: \beta_i=\beta_i^0$ for $i=1,\ldots, r$ and $H_0: \beta_1=\ldots=\beta_r$,
we show that $2[\ell(\widehat{\bs{\beta}}) - \ell(\widehat{\bs{\beta}}^0)]$ converges in distribution to a Chi-square distribution with
the respective degrees of freedoms, $r$ and $r-1$, as the number of nodes $n$ goes to infinity. Here,
$\ell( \bs{\beta})$ is the log-likelihood function on the parameter $\bs{\beta}$, $\widehat{\bs{\beta}}$ is the MLE under the full parameter space,
and $\widehat{\bs{\beta}}^0$ is the restricted MLE under the null parameter space.

\item
For two increasing dimensional null hypotheses $H_0: \beta_i = \beta_i^0$ for $i=1, \ldots, n$ and $H_0: \beta_1=\ldots=\beta_r$ with $r/n \ge c$,
we show that the normalized log-likelihood ratio statistics, $(2[\ell(\widehat{\bs{\beta}}) - \ell(\bs{\beta}^0)] -n)/(2n)^{1/2}$ and
$(2[\ell(\widehat{\bs{\beta}}) - \ell(\widehat{\bs{\beta}}^0)] -r)/(2r)^{1/2}$, both
converge in distribution to the standard normal distribution.
\end{itemize}

Our mathematical arguments depend on the asymptotic expansion of
the log-likelihood function, up to the third order or the fourth order. One key expansion term is involved with a  weighted quadratic sum $\sum_i c_i (d_i - \E(d_i))^2$,
where $d_i$ is the degree of $i$. Note that $d_i$ is not independent over $i$ and classical central limit theorems can not be directly applied.
We prove its central limit theorem by constructing a dependency graph.
Previous results in \cite{simons-yao1999} and \cite{Yan:Xu:2013} including the approximation inverse matrix of the Fisher information matrix and the $\ell_\infty$-error for the MLE are brought here to bound various errors in the remainder terms.

\subsection{Related work}

The $\beta$-model, named by \cite{Chatterjee:Diaconis:Sly:2011},
is an undirected version of the $p_1$ model [\cite{Holland:Leinhardt:1981}] for directed graphs.
The $\beta$-model appears previously in \cite{Park:Newman:2004} and \cite{Blitzstein:Diaconis:2011}.
Since the number of parameters increases with the size of networks,
asymptotic inference is nonstandard.
When the number of nodes goes to infinity,
\cite{Chatterjee:Diaconis:Sly:2011} established the consistency of the MLE in the $\beta$-model;
\cite{Yan:Xu:2013} obtained its asymptotic normality.
\cite{Rinaldo2013} derived the necessary and sufficient conditions of its existence.
\cite{Perry:Wolfe:2012} obtained the finite-sample properties of the maximum likelihood estimate
for a class of $n$-parameter network models with the $\beta$-model as one special case.
The $\beta$-model has been generalized to admit weighted edges [\cite{Hillar:Wibisono:2013,Yan:Zhao:Qin:2015,Yan:Qin:Wang:2015}],
to incorporate a  sparse parameter [\cite{mukherjee2018,Chen:2020}],
to allow covariates [\cite{Wahlstrom:2017,Graham2017}].

For a sparse $\beta$-model that assumes the
connection probability  between nodes $i$ and $j$ takes the form $(\lambda/n)\mu(\beta_i + \beta_j)$ with a common known parameter
$\lambda \in (1, n)$,
\cite{mukherjee2018}  considered the following null hypothesis and alternative hypothesis:
\[
H_0: \bs{\beta} =\mathbf{0} \quad \textrm{vs.} \quad H_1: \bs{\beta} \in \Xi(s_n, A_n) \subset \mathbb{R}_+^n\setminus \{\mathbf{0}\} \label{eqn:hypo}
\]
where $\Xi(s_n,A_n):=\{\beta\in \mathbb{R}_{+}^n: |S(\bs{\beta})|=s_n, \beta_i\ge A, i\in S(\beta)\}$, $S(\bs{\beta}):=\{1\le i\le n:\beta_i\ne 0\}$, and $s_n=n^{1-\alpha}$ with $\alpha \in (0,1)$.
For testing $H_0$, \cite{mukherjee2018} proposed three test statistics: $\sum_i d_i$, $\max_i d_i$ and a higher criticism test based on
$(d_i-\lambda/2)/(\lambda(1-\lambda/2n))^{1/2}$,
and established the asymptotic properties of these statistics under some conditions.
Their problem settings are different from ours.
First, their alternative hypothesis are restricted and only allowed alternative parameters to take positive values
while ours lie in $\R^n$ except for the null parameter space.
Second, their test statistics are not involved with MLEs. Neglecting these information could lead to power loss.
Third, their test statistics requires a known sparse parameter. In practice, this needs to be estimated.
It is unknown whether their results still hold if some estimate $\hat{\lambda}$ is replaced with $\lambda$.

The studies of the Bradley--Terry model in the high dimension setting have attracted great interests in recent years.
In a pioneering work, \cite{simons-yao1999} obtained the upper bound of the $\ell_\infty$-error between the MLE and its true value that directly leads to
the uniform consistency of the MLE, and
established its asymptotic normal distribution under the asymptotic setting where the number of subjects goes to infinity and each pair has the fixed number of comparisons.
\cite{Han-chen2020} extended their results to an Erd\H{o}s--R\'{e}nyi comparison graph under a weak sparsity condition.
The Bradley--Terry model has been widely used for theoretical analysis in various ranking algorithms in machine learning literature
[e.g., \cite{chen2015spectral, negahban2017rank, agarwal2018accelerated,hendrickx2019graph,chen2020partial}],
in which  error bounds for the output estimator $\widehat{\bs{\beta}}$ of the merit parameter $\bs{\beta}$
are established under different conditions/assumptions.
In particular, \cite{chen2019} established the $\ell_\infty$-error for the estimator obtained from the spectral algorithm and the regularized MLE with a $\ell_2$-penalty function.

We note that the $\beta$-model and the Bradley--Terry model can be recast into a logistic regression form.
Although the ``large $N$, diverging $p_N$" framework in generalized linear models (GLMs) has been received much attention in the literature,
most of them focus on the parameter estimation including $M$-estimators [e.g., \cite{huber1973,portnoy1985,BAI1994211,HE2000120}]
and the MLE [e.g., \cite{portnoy1988,haberman1977,ZHOU2021107154}] and generalized estimating equations estimators [\cite{wang2011}].
Here, $N$ is the sample size and $p_N$ is the dimension of parameter space.
Little attention is paid to high dimensional hypothesis testing problems.
\cite{wang2011} obtained the Wilks type of result for the Wald test when $p_N^3/N \to 0$.
In our case, $p_N^3/N \to \infty$, not $0$, where  the dimension of parameter space
 is $p_N = n$  and the total number of observations is $N=n(n+1)/2$ if each edge has only one observation.
In a different setting, by assuming that a sequence of independent and identical distributed samples $\{ X_i\}_{i=1}^n$ from
a regular exponential family $c(\theta) \exp(\theta^\top x)$ with an increasing dimension $p_n$,
\cite{portnoy1988} showed a high dimensional Wilks type of result for the log-likelihood ratio statistic under a simple null hypothesis when $p_N^{3/2}/N \to 0$. Here, our observations are only one dimension and have different distributions.
For logistic regression models with asymptotic regime $p_N/N\in(0, 1/2)$, \cite{sur2019the} showed that the log-likelihood ratio statistic
for testing a single parameter under the null $\beta_i=0$,
converges to a rescaled Chi-square with an inflated factor greater than one.

The rest of the paper is organized as follows.
The Wilks type of theorems for the $\beta$-model  and the Bradley--Terry model are presented in Sections \ref{section-beta-model}
and \ref{section-bt-model}, respectively.
Simulation studies and an application to a NBA data are given in Section \ref{section-numerical}.
All proofs of supported lemmas are relegated to the Supplemental Material.

\input{main.tex}
\input{proof.tex}

\noindent {\bf Acknowledgment}
The views expressed are those of the authors and should not be construed to represent the positions of the Department of the Army or Department of Defense.  J.Z. is partially supported by a National Science Foundation grant and a National Institute of Health grant.

\setlength{\itemsep}{-1.5pt}
\setlength{\bibsep}{0ex}
\bibliographystyle{apa}
\bibliography{reference3}

\end{document}

%% file: main.tex
\section{Wilks theorems for the $\beta$-model}
\label{section-beta-model}

We consider a undirected graph $\mathcal{G}_n$ with $n$ nodes labelled as ``$1, \ldots, n$".
Let $A=(a_{ij})_{n\times n}$ be the adjacency matrix of $\mathcal{G}_n$, where
$a_{ij}$ is an indicator denoting whether node $i$ is connected to node $j$.
That is, $a_{ij}$ is equal to $1$ if there is an edge connecting nodes $i$ and $j$; otherwise, $a_{ij}=0$.
Let $d_i = \sum_{j\neq i} a_{ij}$ be the degree of node $i$ and $\mathbf{d}=(d_1, \ldots, d_n)^\top$ be the degree sequence of $\mathrm{G}_n$.
The $\beta$-model postulates that all $a_{ij}$, $1\le i\neq j \le n$, are mutually independent
Bernoulli random variables, where the probability of $a_{ij}$ being equal to $1$ is $\mu( \beta_i + \beta_j)$ and
$\mu(x)=e^x/(1 + e^x)$ is the logistic function.

The logarithm of the likelihood function under the $\beta$-model can be written as
\begin{equation*}
\ell_\beta(\boldsymbol{\beta}) = \sum_{i,j=1;i\neq j}^n [a_{ij}(\beta_i + \beta_j) - \log(1 + e^{\beta_i + \beta_j})] = \sum_{i=1}^n \beta_i d_i - \sum_{1\le i<j\le n} \log(1 + e^{\beta_i + \beta_j}),
\end{equation*}
where $\boldsymbol{\beta}=(\beta_1, \ldots, \beta_n)$.
As we can see, the $\beta$-model is an undirected exponential random graph model with the degree sequence
as the exclusively natural sufficient statistic.
Setting the derivatives with respect to $\beta_i$ to zero, we obtain the likelihood equations
\begin{equation} \label{eq-likelihood-beta}
d_i = \sum_{j\neq i} \frac{e^{\widehat{\beta}_i + \widehat{\beta}_j}}{1 + e^{\widehat{\beta}_i + \widehat{\beta}_j}},~~i=1,\ldots,n,
\end{equation}
where $\boldsymbol{\widehat{\beta}}=(\widehat{\beta}_1, \ldots, \widehat{\beta}_n)^\top$ is
the MLE of $\boldsymbol{\beta}=(\beta_1, \ldots, \beta_n)^\top$.
The fixed point iterative algorithm in \cite{Chatterjee:Diaconis:Sly:2011} can be used to solve $\boldsymbol{\widehat{\beta}}$.

With some ambiguity of notations, we use $V$ to denote the Hessian matrix of
the negative log-likelihood function under both the $\beta$-model and the Bradley--Terry model.
In the case of the $\beta$-model, the elements of $V$ ($=(v_{ij})_{n\times n}$) are
\begin{equation*}
v_{ii} = \sum_{j\neq i} \frac{e^{\beta_i+\beta_j}}{(1 + e^{\beta_i+\beta_j})^2}, ~~ v_{ij} =  \frac{e^{\beta_i+\beta_j}}{(1 + e^{\beta_i+\beta_j})^2}, ~~i\neq j; ~i,j=1,\ldots, n.
\end{equation*}
$V$ is also the Fisher information matrix of $\bs{\beta}$ and the covariance matrix of $\mathbf{d}$.

We first present Wilks's theorems in the high dimensional setting.
We consider a simple null $H_0: \beta_i=\beta_i^0$, $i=1,\ldots,n$ and a composite null
$H_0: \beta_1=\ldots=\beta_r$, where $\beta_i^0$ for $i=1,\ldots,n$ are  known numbers.
For convenience, we suppress the upper script $0$ in $\beta_i^0$ below.

\begin{theorem}
\label{theorem-LRT-beta}
Define
\[
b_n = \max_{i,j} \frac{(1 + e^{\beta_i+\beta_j})^2}{e^{\beta_i+\beta_j}}.
\]
\begin{itemize}
\item[(a)]
If the following conditions hold:
\begin{equation} \label{condtition-B}
b_n = o\left( \frac{ n^{1/10} }{ (\log n)^{2/5} } \right),~~~~
b_n^3\sum_{i \neq j} |e^{\beta_i + \beta_j} - 1| = o\left( \frac{ n^2 }{ (\log n)^{3/2} } \right),
\end{equation}
then the log-likelihood ratio test statistic $\ell_\beta(\boldsymbol{\hat{\beta}}) - \ell_\beta(\boldsymbol{\beta})$
is asymptotically normally distributed in the sense that
\begin{equation} \label{statistics-beta}
\frac{2\{ \ell_\beta(\boldsymbol{\hat{\beta}}) - \ell_\beta(\boldsymbol{\beta})\} - n}{\sqrt{2n}} \stackrel{d}{\rightarrow} N(0,1), ~~\mbox{as}~~ n\to\infty.
\end{equation}
\item[(b)]
Assume that $r/n\ge \tau>0$, where $\tau$ is a positive constant. If \eqref{condtition-B} holds, then the log-likelihood ratio test statistic
$\ell_\beta(\boldsymbol{\hat{\beta}}) - \ell_\beta(\boldsymbol{\hat{\beta}^{*}})$ is asymptotically normally distributed in the sense that
\begin{equation*} \label{redudant-beta}
\frac{2\{ \ell_\beta(\boldsymbol{\hat{\beta}}) - \ell_\beta(\boldsymbol{\hat{\beta}^{*}})\} - r}{\sqrt{2r}} \stackrel{L}{\rightarrow} N(0,1), ~~\mbox{as}~~ n\to\infty,
\end{equation*}
where $\boldsymbol{\hat{\beta}^{*}}=\arg\max_{\bs{\beta}\in \Theta_1} \ell_\beta(\bs{\beta})$ and $\Theta_1=\{ \bs{\beta}\in \R^n: \beta_1=\cdots=\beta_r\}$.
\end{itemize}
\end{theorem}

We describe briefly the main steps for proving Theorem \ref{theorem-likelihood-ratio-beta} here.
First, we obtain the asymptotic representation of $\widehat{\bs{\beta}}-\bs{\beta}$, which
can be represented as the sum of $V^{-1} (\mathbf{d} - \E \mathbf{d})$ and a remainder term.
Then, we apply a third-order Taylor expansion to  $\ell_\beta(\boldsymbol{\hat{\beta}}$ at point $\bs{\beta}$.
With the use of the maximum likelihood equations,
the first-order and second-order terms can then be expressed as the sum of
$\tfrac{1}{2}(\mathbf{d}-\E \mathbf{d})V^{-1}( \mathbf{d} - \E \mathbf{d} )$ and a remainder term. With the use of the $\ell_\infty$-error for the MLE, the third-order term is asymptotically neglected under condition \eqref{condtition-B}.
The left arguments are to show that $\tfrac{1}{2}(\mathbf{d}-\E \mathbf{d})V^{-1}( \mathbf{d} - \E \mathbf{d} )$ is approximate a Chi-square distribution
with large degrees $n$ of freedom. By using a simple matrix $S=\diag(1/v_{11}, \ldots, 1/v_{nn})$ in \cite{Yan:Xu:2013} to approximate
$V^{-1}$, one can find that its main term is a sum of a sequence of normalized degrees in a weighted quadratic form, i.e.,
$\sum_i (d_i - \E d_i)^2/v_{ii}$. For single $i$, $(d_i - \E d_i)^2/v_{ii}$ is asymptotically a Chi-square distribution because
$d_i$ is a sum of $n-1$ independent Bernoulli random variables. But for all $i$, the terms in the sum are not independent.
Classical central limit theorems can not be directly used.
The central limit theorem for the weighted quadratic sum $\sum_i c_i (d_i - \E(d_i))^2$ is stated below,
which is proved by constructing a dependency graph.

\begin{lemma}\label{lemma:weighte-degree-al}
Consider a general symmetric matrix or a asymmetric matrix $A=(a_{ij})_{n\times n}$ with $|a_{ij}|\le K$ for some constant $K$.
If $A$ is asymmetric, then $a_{ij}+a_{ji}$ is a constant for all $1\le i<j\le n$.
Assume that $c_{\min} \le c_i \le c_{\max}$ for all $i=1,\ldots,n$,
$v_* \le \mathrm{var}( a_{ij} ) \le v_{**}$ for all $i\neq j$, and $\min_{i\neq j}u_{i,j}\ge u_*$, where $u_{ij}=\mathrm{var}( (a_{ij}-\E a_{ij})^2)$.
If the following holds:
\begin{eqnarray}
\label{condition:second}
c_{\max}^6 v_{**}/c_{\min}^5v_*^{5/2}=o(n^{1/2}),~~ u_*/v_*=o(n),
\end{eqnarray}
then the weighted sum $\sum_i c_i (d_i - \E d_i)^2$ is asymptotically normally distributed with mean $\sum_i c_i v_{ii}$
and variance $\sum_i c_i^2 (\sum_{j\neq i}u_{ij}+2v_{ii}^2) + 2\sum_{1\le i<j \le n} c_i c_j u_{ij}$.
\end{lemma}

As a corollary, we have the following result.

\begin{corollary} \label{corollary-weighted-d}
Under the $\beta$-model,
if $b_n=o(n^{1/17})$, then
$\sum_{i=1}^n [d_i-E(d_i)]^2 / v_{ii}$ is asymptotically normally distributed with mean $n$ and variance $2n$.
\end{corollary}

Next, we present the Wilks theorem for a fixed dimensional parameter hypothesis testing problems.
We consider a simple null hypothesis, $H_0: \beta_i=\beta_i^0$, $i=1,\ldots,r$ with a fixed $r$ and a composite null hypothesis,
$H_0: \beta_1=\cdots=\beta_r$.
For convenience, we suppress the upper script $0$ in $\beta_i^0$ below.

\begin{theorem}
\label{theorem-LRT-beta-fixed}
Assume that the conditions in Theorem \ref{theorem-LRT-beta} hold.
\begin{itemize}
\item[(a)]
Under the simple null $H_0: \beta_i=\beta_i^0, i=1, \ldots, r$,
the twice log-likelihood ratio statistic $2(\ell_\beta(\widehat{\bs{\beta}}) - \ell_\beta( \widehat{\bs{\beta}}^0) )$
converges to a Chi-square distribution with $r$ degrees of freedom.

\item[(b)]
Under the composite null $H_0: \beta_1 = \cdots = \beta_r$,
the twice log-likelihood ratio statistic $2(\ell_\beta(\widehat{\bs{\beta}}) - \ell_\beta( \widehat{\bs{\beta}}^0) )$
converges to a Chi-square distribution with $r-1$ degrees of freedom.

\end{itemize}
\end{theorem}

The proof of Theorem \label{theorem-likelihood-ratio-beta-fixed} is much more complex than the proof of Theorem \label{theorem-likelihood-ratio-beta}.
In the fixed dimensional case, it requires to bound $\| \widehat{\bs{\beta}} - \widehat{\bs{\beta}}^0 \|_\infty$ and to evaluate the maximum absolute entry-wise
difference between two approximate inverses matrices.

\section{Wilks theorems for the Bradley--Terry model}
\label{section-bt-model}

In the above, we considered a undirected graph. Now we will consider a weighted directed graph $\mathcal{G}_n$.
As mentioned before, the elements of the adjacency matrix $A$ denote the number of times that one node beats another node.
Let $k_{ij}$ be the number of comparisons between subject $i$ and $j$.
For easy exposition, similar to \cite{simons-yao1999}, we assume $k_{ij}=K$ for all $i\neq j$, where $K$ is a fixed positive constant.
Then, $a_{ij}$ is the number of times that $i$ wins $j$ out of a total number of $K$ comparisons.
The Bradley--Terry model postulates that $a_{ij}$, $1\le i<j \le n$, are mutually independent
binomial random variables, i.e., $a_{ij} \sim \mbox{Binomial}(K, p_{ij})$, where $p_{ij}=\mu(\beta_i - \beta_j)$.
It implies that the win-loss probabilities for any two subjects only depend on the difference of their strength parameters.
The bigger strength parameter, the higher the probability subject $i$ having a win over other subjects.
Let $d_i=\sum_{j\neq i} a_{ij}$ be the total number of wins for subject $j$.

Because the probability is invariable by adding a common constant to all strength parameters $\beta_i$, $i=1, \ldots, n$,
we need a restriction for the identifiability of model. Following \cite{simons-yao1999}, we set $\beta_1=0$ as a constraint.
Notice that the number of free parameters here is $n-1$, different from the $\beta$-model with $n$ free parameters.
The logarithm of the likelihood function under the Bradley--Terry model is
\begin{equation}\label{likelihood-bt}
\ell_{bt}(\boldsymbol{\beta}) = \sum_{i,j=1;i\neq j}^n a_{ij} [\beta_i - \log(e^{\beta_i}+e^{\beta_j})] = \sum_{i=1}^n \beta_id_i  - K \sum_{1\le i<j\le n} \log(e^{\beta_i}+e^{\beta_j}),
\end{equation}
where $\boldsymbol{\beta} = (\beta_2, \ldots, \beta_n)^\top$ and $\beta_1=0$.
As we can see, it is a directed exponential random graph model with the out-degree sequence as its natural sufficient statistic.
Setting the derivatives with respect to $\beta_i$ to zero, we obtain the likelihood equations
\begin{equation} \label{estimated-eq-bt-a}
d_i = \sum_{j=1,j\neq i}^n \frac{K e^{\hat{\beta}_i}}{e^{\hat{\beta}_i} + e^{\hat{\beta}_j}}, ~~i=2,\ldots,n,
\end{equation}
where $\boldsymbol{\widehat{\beta}} = (\widehat{\beta}_2, \ldots, \hat{\beta}_n)$ is the maximum likelihood estimate of $\boldsymbol{\beta}$ with $\widehat{\beta}_1=0$.
If the directed graph $\mathcal{G}_n$ is strongly connected, then the MLE exists and is unique [\cite{Ford1957}].
Note that $d_1$ is not involved in \eqref{estimated-eq-bt-a}; indeed, given $d_2, \ldots, d_n$ and $K$, $d_1$ is determined.

Now, we present the Wilks type of theorems for the Rradley--Terry model.
The proof is similar to the proof of Theorem \ref{theorem-likelihood-ratio-beta} and is put in the Supplementary Material.

\begin{theorem} \label{theorem-likelihood-ratio-bt}
(a) If the following conditions hold:
\begin{equation} \label{bt-condition1}
M_n = o ( n^{1/17} ),~~~~
\sum_{i,j=1}^n \left| \frac{e^{\beta_i} - e^{\beta_j}}{e^{\beta_i} + e^{\beta_j}} \right| = o \left( \frac{ n^{25/14} }{ (\log n)^{15/7} } \right),
\end{equation}
then the log-likelihood ratio test statistic $\ell_{bt}(\boldsymbol{\hat{\beta}}_n) - \ell_{bt}(\boldsymbol{\beta}_n)$ is asymptotically normally distributed in the sense that
\begin{equation} \label{statistics-bt1}
\frac{2\{ \ell_{bt}(\boldsymbol{\hat{\beta}}_n) - \ell_{bt}(\boldsymbol{\beta}_n) \} - (n-1)}{\sqrt{2(n-1)}}
\stackrel{L}{\rightarrow} N(0,1), ~~\mbox{as}~~ n\to\infty.
\end{equation}
(b)Without loss of generality, suppose the composite null hypothesis takes the following form, i.e. $H_0^*$: $\beta_1 = \cdots = \beta_r$, $2 \leq r \leq n$.  Let $\boldsymbol{\hat{\beta}}_{n}^* = (\hat{\beta}_{1}^{*}, \hat{\beta}_{2}^{*}, \ldots, \hat{\beta}_{n}^{*})$ be the maximum
likelihood estimate of $\boldsymbol{\beta}_n$ under $H_0^*$, with $\hat{\beta}_1^*=0$. Assume that $r/n \ge \tau>0$, where $\tau$ is a positive constant. If \eqref{bt-condition1} holds, then the log-likelihood ratio test statistic
$\ell_{bt}(\boldsymbol{\hat{\beta}}_n) - \ell_{bt}(\boldsymbol{\hat{\beta}}_n^{*})$ is asymptotically normally distributed in the sense that
\begin{equation} \label{redudant-bt}
\frac{2\{ \ell_{bt}(\boldsymbol{\hat{\beta}}_n) - \ell_{bt}(\boldsymbol{\hat{\beta}}_{n}^*)\} - (r-1)}{\sqrt{2(r-1)}} \stackrel{L}{\rightarrow} N(0,1), ~~\mbox{as}~~ n\to\infty.
\end{equation}
\end{theorem}

Note that in order to guarantee the existence of the maximum likelihood estimate with high probability,
it is necessary to control the increasing rate of $M_n$ as discussed in \cite{simons-yao1999}.
In the case that some $u_i$'s are
very large while others are very small, corresponding to a large value of $M_n$,
the nodes with relatively poor merits will stand very little chance of beating those
with relatively large merits.
Whenever all nodes could be partitioned into two sets,
in which the vertices in one set will win all games against those in the other set,
the MLE will not exist [\cite{Ford1957}].
Therefore, the first condition in \eqref{bt-condition1} is used to control the increasing rate of $M_n$.
Moreover, it would be of interest to see whether the condition imposed on $M_n$ can be relaxed.
The second condition in \eqref{bt-condition1} is technical, due to the control of the remainder in the Taylor expansion of the log-likelihood function, which essentially requires that the merits of different vertices do not differ too much.
$\tau$ is used to control the number of parameters that are equal. A larger $\tau$, more parameters are equal.

Note that in the above discussion, we have assumed the $k_{ij}$'s, $i\neq j$ are all equal to a constant $K$.  This is only for the purpose of simplifying notations. Theorem \ref{likelihood-ratio-bt} can be readily extended to the general case,
where $k_{ij}$'s are not necessarily the same (but with a bound).
A complicated case is when the $k_{ij}$
are quite different from one another. For example, a large number of pairs don't have direct comparisons
or some pairs have too many comparisons. In these cases, it is unknown whether Theorem \ref{likelihood-ratio-bt}
still holds. We would like to investigate this problem in the future work.

Next, we present the Wilks theorem for a fixed dimensional parameter hypothesis testing problems.

\begin{theorem}
\label{theorem-likelihood-ratio-bt-fixed}
Assume that the conditions in Theorem \ref{theorem-likelihood-ratio-bt} hold.
\begin{itemize}
\item[(a)]
Under the simple null $H_0: \beta_i=\beta_i^0, i=1, \ldots, r$,
the twice log-likelihood ratio statistic $2(\ell_{bt}(\widehat{\bs{\beta}}) - \ell_{bt}( \widehat{\bs{\beta}}^0) )$
converges to a Chi-square distribution with $r-1$ degrees of freedom.

\item[(b)]
Under the composite null $H_0: \beta_1 = \cdots = \beta_r$,
the twice log-likelihood ratio statistic $2(\ell_{bt}(\widehat{\bs{\beta}}) - \ell_{bt}( \widehat{\bs{\beta}}^0) )$
converges to a Chi-square distribution with $r-2$ degrees of freedom.

\end{itemize}
\end{theorem}

\section{Numerical Results}
\label{section-numerical}
In this section, we demonstrate the theoretical results via numerical studies.

\subsection{Simulation studies}
We carry out simulations to evaluate the performance of the log-likelihood ratio statistics for finite number of nodes.
We only present the simulation results under the Bradley--Terry model here and those for the $\beta$-model are in the supplementary material.
To evaluate Theorems \ref{theorem-likelihood-ratio-bt}, we considered several simulations.
In all simulation studies, we let the number of experiments $k_{ij}$ equal to 1 for all $1\le i\neq j\le n$, and the parameters $\beta_i$, $i=1,\ldots, n$, be in a linear form.  Specifically, for the simple null, we set $\beta_i=(i-1)L_n/(n-1)$, $i=1,\ldots,n$, and for the composite null, we set $\beta_1=\cdots=\beta_r=0$, where $r=n/2$ and $\beta_i=(i-1)L_n/(n-1)$, $i=r+1, \ldots, n$.  Note that in both settings of $\beta_i$'s, $L_n=\max_i \beta_i$ and $M_n=e^{L_n}$.
Four values of $L_n$ were considered, specifically, 0, $\log(\log n)$, $\log n$ and $n$, and consequently $M_n=1$, $\log n$, $n$ and $e^n$ respectively.  In each simulation, we computed the test statistic as described in the corresponding theorem, and the procedure was repeated $10,000$ times.

The results for the simple null and composite null in the Bradley-Terry model are shown in Figures 1 and 2,
respectively.
In each QQ-plot, the horizontal and vertical axes correspond to the theoretical and empirical quantiles respectively.  Note that when $M_n=e^n$, condition \eqref{bt-condition1} is not satisfied, and we observed that the maximum likelihood estimate did not exist more than $90$\% times out of the $10,000$ repetitions, thus the corresponding result is not reported; on the other hand, the maximum likelihood estimate always existed for other values of $M_n$, i.e. 1, $\log n$ and $n$,
which is in agreement with earlier findings in \cite{simons-yao1999}.
As we can see, when $n=50$, the empirical quantiles differ a little from the theoretical ones, but as $n$ increases to $500$, the difference diminishes and  the empirical quantiles agree well with the theoretical ones.  Further,
we can see that as $M_n$ increases, the difference between the empirical quantiles and the theoretical ones becomes more prominent.

\begin{figure}[!htb]
\centering
\caption{Simulation results for the Bradley-Terry model under the simple null.  The horizontal and vertical axes in each QQ-plot are the theoretical (based on the standard normal distribution) and empirical quantiles (based on the log-likelihood ratio test statistic), respectively.  The straight lines correspond to $y=x$.
The first, second, third columns correspond to $M_n=1, \log n, n$, respectively.
}
\subfigure[$n=50$]{ \raisebox{-1cm}{\includegraphics[width=1.0\textwidth, height=4cm]{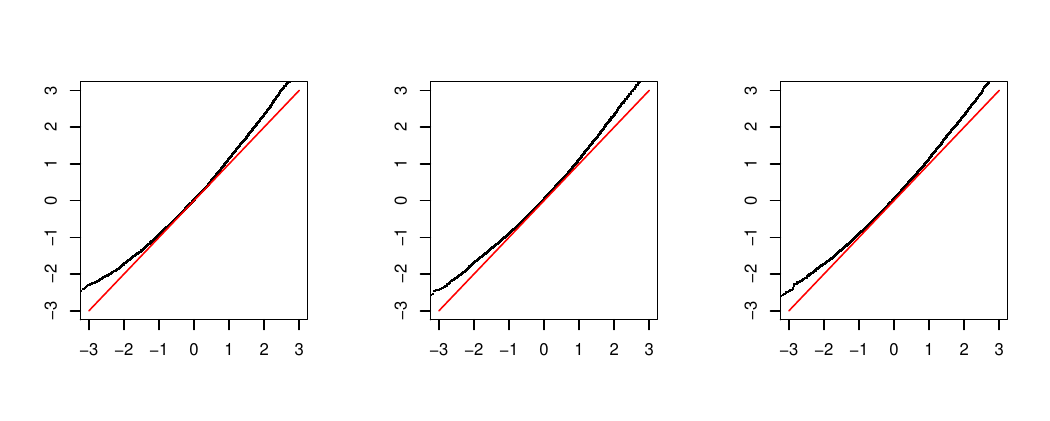}}}
\subfigure[$n=200$]{ \raisebox{-1cm}{\includegraphics[width=1.0\textwidth, height=4cm]{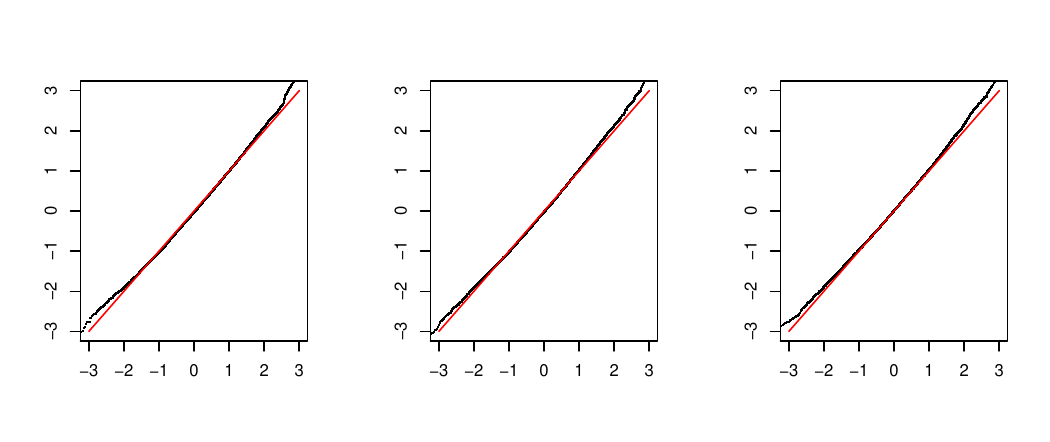} }}
\subfigure[$n=500$]{ \raisebox{-1cm}{\includegraphics[width=1.0\textwidth, height=4cm]{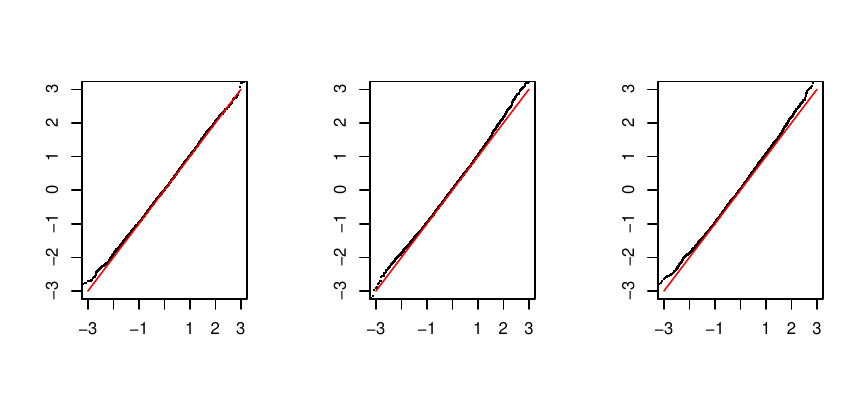} }}
\label{fig:Bradley-Terry}     
\end{figure}

\begin{figure}[!htb]
\centering
\caption{Simulation results for the Bradley-Terry model under the composite null.  The horizontal and vertical axes in each QQ-plot are the theoretical (based on the standard normal distribution) and empirical quantiles (based on the log-likelihood ratio test statistic), respectively.  The straight lines correspond to $y=x$.
The first, second, third columns correspond to $M_n=1, \log n, n$, respectively.
}
\subfigure[$n=50$, $r=25$]{ \raisebox{-1cm}{\includegraphics[width=1.0\textwidth, height=4cm]{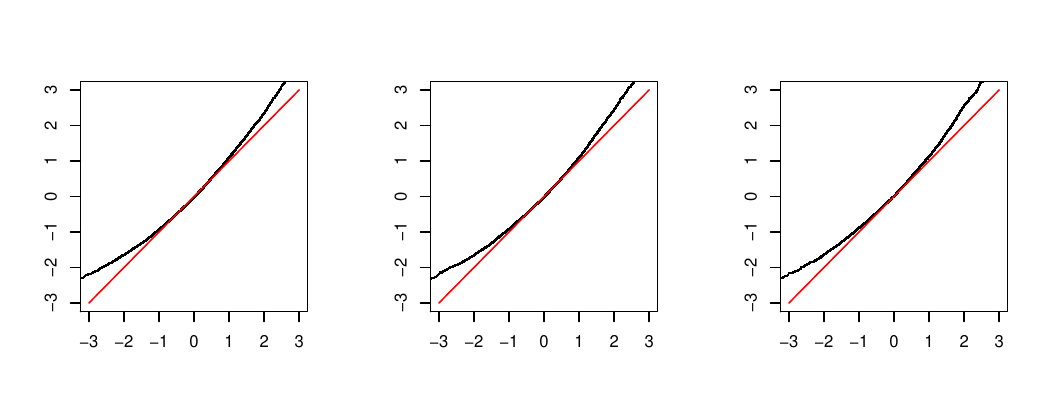} }}
\subfigure[$n=200$, $r=100$]{ \raisebox{-1cm}{\includegraphics[width=1.0\textwidth, height=4cm]{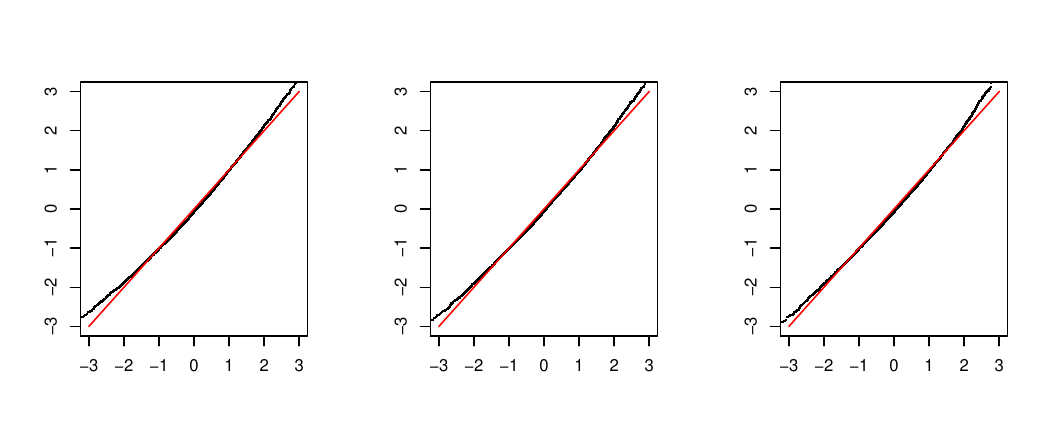} } }
\subfigure[$n=500$, $r=250$]{ \raisebox{-1cm}{\includegraphics[width=1.0\textwidth, height=4cm]{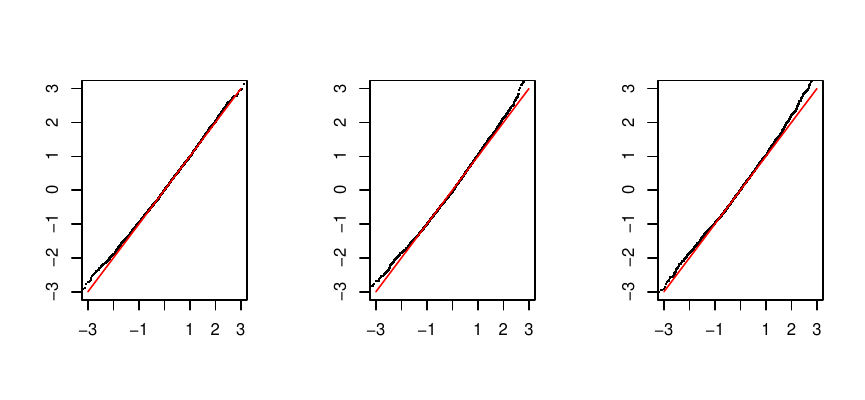} } }
\label{fig:Bradley-Terry}     
\end{figure}

Next, we investigate the powers of the test statistic \eqref{table-bt-powers}.
The null took the form $H_0: \beta_1=\cdots=\beta_r=0$, and the true model was set to be $\beta_i=ic/r$, $i=1,\ldots,r$.  The other parameters were set as $\beta_i=(i-r)L_n/n$ for $i=r+1,\ldots,n$.  The results are shown in Table \ref{table-bt-powers}.  We can see that when $c=0$, the simulated type I errors agree reasonably well with the nominal level, even when $n=50$.  Further, when $n$ and $r$ are fixed, as $c$ increases, the power tends to increase.  Similar phenomenon can be observed when $r$ increases while $n$ and $c$ are fixed, or when $n$ increases while $c$ and $r$ are fixed.

\begin{table}[h!]
\centering
\caption{Powers of the proposed likelihood ratio tests} \label{table-bt-powers} \vskip5pt
\small
\begin{tabular}{ccc ccc cc}
\hline
\\
\multicolumn{8}{c}{Powers of the test (5) for the Bradley-Terry model}\\ 
\hline
\\
$n$    & $L_n$        & $r$  &$c=0$       & $c=0.4$ & $c=0.8$ & $c=1.2$ & $c=1.6$ \\
\hline
$n=30$ & $0$          & $10$ &$5.66$      &$ 9.53$&$ 23.28$&$ 48.40 ( 3\times 10^{-4} )$&$ 74.61 ( 1.2\times 10^{-3} )$                 \\
       &              & $20$ &$5.88 ( 0 )$&$ 10.26 ( 0 )$&$ 32.14 ( 0 )$&$ 69.23 ( 0 )$&$ 93.92 ( 2\times 10^{-4} )$         \\
       &$\log(\log n)$& $10$ &$5.66 ( 0 )$&$ 9.27 ( 0 ) $&$ 23.67 ( 0 )$&$ 50.55 ( 0 )$&$ 78.24 ( 1\times 10^{-4} )$         \\
       &              & $20$ &$5.69 ( 0 )$&$ 10.38 ( 0 )$&$ 32.51 ( 0 )$&$ 69.56 ( 0 )$&$ 94.45 ( 0 )$             \\
$n=50$ & $0$          & $10$ &$5.09 ( 0 )$&$ 11.07 ( 0 )$&$ 36.19 ( 0 )$&$ 71.96 ( 0 )$& 93.67 $( 1\times 10^{-4} )$         \\
       &              & $20$ &$5.18 ( 0 )$&$ 13.86 ( 0 )$&$ 53.32 ( 0 )$&$ 92.22 ( 0 )$&$ 99.75 ( 0 )$             \\
       &$\log(\log n)$& $10$ &$5.21 ( 0 )$&$ 11.07 ( 0 )$&$ 36.88 ( 0 )$&$ 75.76 ( 0 )$&$ 96.06 ( 0 )$\\
       &              & $20$ &$5.23 ( 0 )$&$ 13.32 ( 0 )$&$ 54.58 ( 0 )$&$ 93.70 ( 0 )$&$ 99.86 ( 0 )$ \\
     \hline
\end{tabular}
\end{table}

\subsection{A data example}
National Basketball Association (NBA) is one of the most successful men's professional basketball league in the world.
The current league organization divides its total thirty teams into two conferences: the western conference and the eastern conference.
In the regular season, every team plays with every other team three or four times.
It would be of interest to test whether there are significant difference among a set of teams.
Here we use the 2008-09 NBA season data as an illustrative example.

The fitted merits using the Bradley--Terry model are presented in Table \ref{table-merits}, in which Philadelphia 76ers is the reference team.
As we can see, the ranking based on the won-loss percentage and that based on the fitted merits are similar.
Further, we use \eqref{redudant-bt} to test whether there are significant differences among the middle 9 teams according to the ranking of the
won-loss percentage, i.e. No. 4--12, in each conference.
It may be obvious that it is significance if testing the equality of all teams in each conference.
In fact, we get asymptotic p-values in the magnitude of $10^{-86}$. So we drop the top $3$ teams and bottom $3$ ones.
The values of \eqref{redudant-bt} are $0.290$ and $13.6$ for the eastern conference and the western conference respectively, with the corresponding p-values $0.772$ and $2.8\times 10^{-42}$.
To evaluate the quality of asymptotic approximation, we used the permutation tests under the null based on $100,000$ Monte Carlo simulations, getting the p-values $0.739$ and $<10^{-5}$. We can see that the empirical one and the asymptotic one are similar for
testing the equality of the middle $9$ teams in the east conference. For the west conference, the asymptotic one gives much smaller p-value.
The results indicate that there is no significant difference among the middle nine teams in the eastern conference while there are significant differences among the those teams in the western conference.

\begin{table}[h!]
\centering
\caption{Fitted merits based on the 2008-09 NBA season data. The values in parentheses are the standard errors.} \label{table-merits} \vskip5pt
\begin{tabular}{l llc | llc }
\hline
& \multicolumn{3}{c|}{  Eastern Conference } & \multicolumn{3}{c}{  Western Conference  }\\
&    Team      &  W-L    & Merit          &     Team      & W-L       & Merit   \\
\hline
1&  Cleveland Cavaliers  &  66-16   &4.532(0.374)  & Los Angeles Lakers      &65-17 &  4.158(0.370) \\
2&  Boston Celtics       &  62-20   &3.462(0.359)   & Denver Nuggets          &54-28 &  2.058(0.344) \\
3&  Orlando Magic        &  59-23   &2.745(0.351)   & San Antonio Spurs       &54-28 &  2.005(0.343) \\
4&  Atlanta Hawks        &  47-35   &1.404(0.337)   & Portland Trail Blazers  &54-28 &  2.059(0.344) \\
5&  Miami Heat           &  43-39   &1.146(0.335)   & Houston Rockets         &53-29 &  1.953(0.342) \\
6&  Philadelphia 76ers   &  41-41   &1.000   & Dallas Mavericks        &50-32 &  1.612(0.339) \\
7&  Chicago Bulls        &  41-41   &1.002(0.334)   & New Orleans Hornets     &49-33 &  1.563(0.339) \\
8&  Detroit Pistons      &  39-43   &0.899(0.335)   & Utah Jazz               &48-34 &  1.425(0.338) \\
9&  Indiana Pacers       &  36-46   &0.794(0.335)   & Phoenix Suns            &46-36 &  1.284(0.338) \\
10& Charlotte Bobcats    &  35-47   &0.716(0.335)   & Golden State Warriors   &29-53 &  0.502(0.343) \\
11& New Jersey Nets      &  34-48   &0.682(0.336)   & Minnesota Timberwolves  &24-58 &  0.383(0.351) \\
12& Milwaukee Bucks      &  34-48   &0.697(0.336)   & Memphis Grizzlies       &24-58 &  0.387(0.387) \\
13& Toronto Raptors      &  33-49   &0.659(0.337)   & Oklahoma City Thunder   &23-59 &  0.349(0.353) \\
14& New York Knicks      &  32-50   &0.621(0.338)   & Los Angeles Clippers    &19-63 &  0.272(0.364) \\
15& Washington Wizards   &  19-63   &0.283(0.361)   & Sacramento Kings        &17-65 &  0.230(0.371) \\
\hline
\end{tabular}

\end{table}


%% file: proof.tex
\section{Proofs}
We introduce some notations.
For a vector $\mathbf{x}=(x_1, \ldots, x_n)^\top\in \R^n$, denote by $\|\mathbf{x}\|$ for a general norm on vectors with the special cases
$\|\mathbf{x}\|_\infty = \max_{1\le i\le n} |x_i|$ and $\|\mathbf{x}\|_1=\sum_i |x_i|$ for the $\ell_\infty$- and $\ell_1$-norm of $\mathbf{x}$ respectively.
For an $n\times n$ matrix $J=(J_{ij})$, let $\|J\|_\infty$ denote the matrix norm induced by the $\ell_\infty$-norm on vectors in $\R^n$, i.e.,
\[
\|J\|_\infty = \max_{\mathbf{x}\neq 0} \frac{ \|J\mathbf{x}\|_\infty }{\|\mathbf{x}\|_\infty}
=\max_{1\le i\le n}\sum_{j=1}^n |J_{ij}|,
\]
and $\|J\|$ be a general matrix norm.

We define a matrix class $\mathcal{L}_n(b_{n0}, b_{n1})$ with two positive numbers $b_{n0}$ and $b_{n1}$.
We say an $n\times n$ matrix $V=(v_{ij})$ belongs to the matrix class
$\mathcal{L}_n(b_{n0}, b_{n1})$ if
\[
\begin{array}{cl}
v_{ii}=\sum_{j\neq i} v_{ij}, & i=1, \ldots, n \\
b_{n0} \le v_{ij} \le b_{n1}, & i,j=1, \ldots, n; i\neq j.
\end{array}
\]
For $V\in \mathcal{L}_n(b_{n0}, b_{n1})$, \cite{hillar2012inverses} obtained a tight bound of $\|V\|_\infty$.
As applied here, we have that for $V\in \mathcal{L}_n(1/b_n, 1/4)$ and $n\ge 3$,
\begin{equation}\label{ineq-tight-V}
\frac{2}{(n-1)} \le \|V^{-1}\|_\infty \le \frac{(3n-4)b_n}{2(n-1)(n-2)}.
\end{equation}
\cite{Yan:Xu:2013} proposed to use a simple matrix $\mathrm{diag}(1/v_{11}, \ldots, 1/v_{nn})+1/v_{\cdot\cdot}$ to approximate $V^{-1}$, where
$v_{\cdot\cdot}=\sum_i v_{ii}$. Since $v_{ii}/v_{\cdot\cdot} =  O(1/n)$, we use the diagonal matrix
\begin{equation}\label{definition-S}
S=\mathrm{diag}(1/v_{11}, \ldots, 1/v_{nn}),
\end{equation}
as an approximation for further simplification. They proved
\begin{equation}\label{ineq-V-S-appro-upper-b}
\|V^{-1} - S \|_{\max} =O\left(\frac{b_n^3}{n^2}\right).
\end{equation}

Recall that $\mu(x) =  e^x/(1 + e^x)$.
A direct calculation gives that the derivative of $\mu(x)$ up to the third order are
\begin{eqnarray}\label{eq-derivative-mu-various}
\mu^\prime(x) = \frac{e^x}{ (1+e^x)^2 },~~  \mu^{\prime\prime}(x) = \frac{e^x(1-e^x)}{ (1+e^x)^3 },~~ \mu^{\prime\prime\prime}(x) =  \frac{ e^x [ (1-e^x)^2 - 2e^x] }{ (1 + e^x)^4 }.
\end{eqnarray}
It is not difficult to verify the following inequalities:
\begin{equation}\label{ineq-mu-deriv-bound}
|\mu^\prime(x)| \le \frac{1}{4}, ~~ |\mu^{\prime\prime}(x)| \le \frac{1}{4},~~ |\mu^{\prime\prime\prime}(x)| \le \frac{1}{4}.
\end{equation}
This fact will be used in the proofs repeatedly.
Let $\bar{a}_{ij} = a_{ij} - \E(a_{ij})$ be the centered random variable of $a_{ij}$ and define $\bar{a}_{ii}=0$ for all $i=1, \ldots, n$.
Correspondingly, denote $\bar{d}_i = d_i - \E(d_i)$ and $\bar{\dd}=(\bar{d}_1, \ldots, \bar{d}_n)^\top$.

\subsection{Proofs under the $\beta$-model}
\label{section:beta}

Let $V = -\partial^2 \ell_{\beta}( \bs{\beta} )/\partial \bs{\beta} \partial \bs{\beta}^\top$.
It is easy to see that the elements $v_{ij}$ of $V$, $i,j=1,\ldots, n$, are
\begin{equation*}
v_{ii} = \sum\nolimits_{j\neq i} v_{ij}, ~~ v_{ij} = \frac{e^{\beta_i + \beta_j}}{(1 + e^{\beta_i + \beta_j})^2}=\mu^\prime(\beta_i + \beta_j), ~~i \neq j; i,j=1,\ldots, n,
\end{equation*}
which is also the covariance matrix of $\mathbf{d}$ and the Fisher information matrix of $\bs{\beta}$.
To prove Theorem \ref{theorem-LRT-beta}, we need the upper bound of $\bar{\dd}$ and
the error bound $\|\widehat{\bs{\beta}}-\bs{\beta}^0\|_\infty$, stated as two lemmas below.
In this section, we will suppress the superscript $0$ in $\bs{\beta}^0$ for convenience when causing no confusions.

\begin{lemma} \label{lemma-clt-beta-W}
Let $W=V^{-1}-S$.
If $b_n=o( n^{1/6} )$, then $\bar{\dd}^\top W \bar{\dd} = o_p(n)$.
\end{lemma}

\begin{lemma}\label{lemma-consi-beta}
If $b_n = o( (n/\log n)^{1/4})$, then with probability at least $1-O(n^{-1})$, the MLE $\widehat{\bs{\beta}}$ exists and satisfies
\[
\| \widehat{\bs{\beta}} - \bs{\beta} \|_\infty  = O\left( b_n\sqrt{ \frac{\log n}{n} } \right).
\]
\end{lemma}

We are now ready to prove the first part of Theorem \ref{theorem-LRT-beta}.

\begin{proof}[Proof of Theorem \ref{theorem-LRT-beta} (a)]
Let $E_n$ be the event that the MLE in \eqref{eq-likelihood-beta} exists and satisfies that
\begin{equation}\label{bound-expij-beta}
\| \bs{\widehat{\beta}} - \bs{\beta} \|_\infty  \lesssim \left( b_n\sqrt{ \frac{\log n}{n} } \right).
\end{equation}
By Lemma \ref{lemma-consi-beta}, the event $E_n$ holds with probability at least $1-O(n^{-1})$ if $b_n=o((n/\log n)^{1/4})$.
The following calculations are based on the event $E_n$.
The proof of Theorem \ref{theorem-likelihood-ratio-beta} (a) is divided into three steps.
Step 1 is about the asymptotic representation of $\widehat{\bs{\beta}}$.
Step 2 is about the asymptotic expansion of $\ell_{\beta}( \bs{\widehat{\beta}} )$.
Step 3 is a combination step.

Step 1. We characterize the asymptotic representation of $\widehat{\bs{\beta}}$.
To simplify notations, define
 $\widehat{\pi}_{ij} = \widehat{\beta}_i + \widehat{\beta}_j$ and $\pi_{ij}= \beta_i + \beta_j$.
A second-order Taylor expansion gives that
\begin{eqnarray}\label{eq-expansion-beta-a}
\mu( \widehat{\pi}_{ij} )
&=& \mu( \pi_{ij} ) + \mu^\prime(\pi_{ij}) (\widehat{\pi}_{ij} - \pi_{ij}) +
\frac{1}{2} \mu^{\prime\prime}( \tilde{\pi}_{ij} ) (\widehat{\pi}_{ij} - \pi_{ij})^2,
\end{eqnarray}
where $\tilde{\pi}_{ij}$ lies between $\widehat{\pi}_{ij}$ and $\pi_{ij}$.
Let
\begin{equation}\label{eq:definition:h}
h_{ij}= \frac{1}{2}\mu^{\prime\prime}( \tilde{\pi}_{ij} )(\widehat{\pi}_{ij} - \pi_{ij})^2,~~
h_i=\sum_{j\neq i}h_{ij}, ~~\mathbf{h}=(h_1, \ldots, h_n)^\top.
\end{equation}
In view of \eqref{ineq-mu-deriv-bound} and \eqref{bound-expij-beta}, we have
\begin{equation}\label{ineq-beta-h}
\| \mathbf{h} \|_\infty \le \frac{1}{2}(n-1) \max_{i,j} |h_{ij} | < n\| \bs{\widehat{\beta}} - \bs{\beta} \|_\infty^2
\lesssim b_n^2 \log n.
\end{equation}
By \eqref{eq-likelihood-beta} and \eqref{eq-expansion-beta-a}, we have
\begin{equation*}
d_i-\E(d_i)=\sum_{j=1,j\neq i}^n v_{ij}\{ (\widehat{\beta}_i-\beta_i)+(\widehat{\beta}_j-\beta_j)\} + h_i, ~~~i=1, \ldots, n.
\end{equation*}
Writing the above equations into the matrix form, we have
\begin{equation*}
\mathbf{d} - \E( \mathbf{d} ) = V ( \widehat{\boldsymbol{\beta}} - \boldsymbol{\beta} ) + \mathbf{h}.
\end{equation*}
Therefore,
\begin{equation}\label{eq-expansion-hatbeta-beta}
\boldsymbol{\widehat{\beta}} - \boldsymbol{\beta} = V^{-1}\bar{\dd} - V^{-1}\mathbf{h},
\end{equation}
where, by \eqref{ineq-tight-V} and \eqref{ineq-beta-h},
\begin{equation}\label{ineq-beta-h-b}
\| V^{-1}\mathbf{h} \|_\infty \le \| V^{-1} \|_\infty \| \mathbf{h} \|_\infty \lesssim \frac{ b_n^3 \log n }{ n }.
\end{equation}

Step 2. We derive the asymptotic expansion of $\ell_{\beta} ( \boldsymbol{\widehat{\beta}})$.
Applying a third-order Taylor expansion to $\ell_{\beta} ( \boldsymbol{\widehat{\beta}})$ gives that
 \renewcommand{\arraystretch}{1.5}
\begin{equation}\label{Taylor-ell-beta-a}
\begin{array}{rcl}
\ell_{\beta}( \bs{\widehat{\beta}} ) & = & \ell_{\beta}( \bs{\beta} )
+ \frac{\partial \ell_{\beta}( \bs{\beta} ) }{\partial \bs{\beta}^\top } ( \bs{\widehat{\beta} } - \bs{\beta})
+ \frac{1}{2} ( \bs{\widehat{\beta} } - \bs{\beta})^\top \frac{\partial^2 \ell_{\beta}( \bs{\beta} )
}{\partial \bs{\beta} \partial \bs{\beta}^\top }( \bs{\widehat{\beta} } - \bs{\beta}) \\
&&+ \frac{1}{6} \sum_{i,j,k} \frac{\partial^3 \ell_{\beta}( \widetilde{\bs{\beta}} )
}{\partial \beta_i\partial \beta_j \partial \beta_k } (\widehat{\beta}_i-\beta_i)(\widehat{\beta}_j-\beta_j)(\widehat{\beta}_k-\beta_k),
\end{array}
\end{equation}
where $\widetilde{\bs{\beta}}=\bs{\beta}+t\bs{\widehat{\beta}}$ for some $t\in (0, 1)$.
For the third-order expansion term, observe that for three distinct indices $i, j, k$, we have
\[
 \frac{\partial^3 \ell_{\beta}( \tilde{\bs{\beta}} )
}{\partial \beta_i\beta_j\beta_k } = 0,
\]
and
\begin{equation*}  
\frac{\partial^3 \ell_{\beta}(\boldsymbol{\beta})}{\partial \beta_i^3 }=\sum_{j\neq i} \mu^{\prime\prime}( \pi_{ij} ),~~~~
\frac{\partial^3 \ell_{\beta}(\boldsymbol{\beta})}{\partial \beta_i^2\partial \beta_j }= \mu^{\prime\prime}( \pi_{ij} ).
\end{equation*}
It follows that
\renewcommand{\arraystretch}{1.5}
\begin{equation}\label{cacluation-z}
\begin{array}{rcl}
z : & = & \sum_{i,j,k} \frac{\partial^3 \ell_{\beta}( \widetilde{\bs{\beta}} )
}{\partial \beta_i\partial \beta_j \partial \beta_k } (\widehat{\beta}_i-\beta_i)(\widehat{\beta}_j-\beta_j)(\widehat{\beta}_k-\beta_k) \\
& = & \frac{1}{6}\{ \sum_{i=1}^n  (\widehat{\beta}_i-\beta_i)^3 \sum_{j\neq i} \mu^{\prime\prime}( \tilde{\pi}_{ij} )
+ 2 \sum_{i,j=1, j\neq i}^n  (\widehat{\beta}_i-\beta_i)^2(\widehat{\beta}_j-\beta_j)\mu^{\prime\prime}( \tilde{\pi}_{ij} ) \}.
\end{array}
\end{equation}
Therefore, the difference between $\ell_{\beta}( \bs{\widehat{\beta}} )$ and $\ell_{\beta}( \bs{\beta} )$ can be expressed as
\begin{equation}\label{lrt-a-beta}
\ell_{\beta}( \bs{\widehat{\beta}} ) - \ell_{\beta}( \bs{\beta} )
= ( \bs{\widehat{\beta} } - \bs{\beta})^\top (\mathbf{d} - \E \mathbf{d})
+ \frac{1}{2} ( \bs{\widehat{\beta} } - \bs{\beta})^\top V( \bs{\widehat{\beta} } - \bs{\beta}) + \frac{1}{6}z.
\end{equation}

Step 3. We combine aforementioned two steps and bound remainder terms.
Substituting \eqref{eq-expansion-hatbeta-beta} into \eqref{lrt-a-beta}, it yields
\begin{equation}\label{minus-ratio-b}
\ell_{\beta}(\boldsymbol{\widehat{\beta}})-\ell_{\beta}(\boldsymbol{\beta})=
\frac{1}{2}\{ \mathbf{d} - \E (\mathbf{d}) \}^\top V^{-1} \{ \mathbf{d} - \E( \mathbf{d} ) \}
 - \frac{1}{2}\mathbf{h}^\top V^{-1} \mathbf{h} + z.
\end{equation}
In view of Corollary \ref{corollary-weighted-d} and Lemma \ref{lemma-clt-beta-W},  the remainder of the proof is to show
\begin{equation}\label{th22}
\frac{ \mathbf{h}^\top V^{-1} \mathbf{h} }{\sqrt{n}}=o_p(1),~~\frac{z}{\sqrt{n}}=o_p(1).
\end{equation}
By \eqref{ineq-beta-h} and \eqref{ineq-beta-h-b}, we have
\[
|\mathbf{h}^\top V^{-1} \mathbf{h}| \le n \| \mathbf{h} \|_\infty \| V^{-1} \mathbf{h}| \|_\infty
\lesssim
 n \cdot b_n^2\log n \cdot \frac{ b_n^3 \log n}{n} \lesssim b_n^5 (\log n)^2.
\]
If $b_n =o( n^{1/10}/(\log n)^{2/5})$, then
\[
\frac{1}{n^{1/2}} |\mathbf{h}^\top V^{-1} \mathbf{h}| \lesssim \frac{ b_n^5(\log n)^2 }{n^{1/2}}  = o(1).
\]
We now bound $z$.
By the mean value theorem, we have
\[
\mu^{\prime\prime}( \tilde{\pi}_{ij} ) = \mu^{\prime\prime}(\pi_{ij})
+ \mu^{\prime\prime\prime}(\bar{\pi}_{ij} )(\tilde{\pi}_{ij} - \pi_{ij}),
\]
where $\bar{\pi}_{ij}$ lies between $\pi_{ij}$ and $\tilde{\pi}_{ij}$.
By \eqref{ineq-mu-deriv-bound},  we have
\[
\sum_{i\neq j } | \mu^{\prime\prime}(\tilde{\beta}_i + \tilde{\beta}_j) |
\le \sum_{i\neq j } | \mu^{\prime\prime}(\beta_i + \beta_j) | + \frac{n(n-1)}{4} \|\bs{\widehat{\beta}}- \bs{\beta} \|_\infty.
\]
It follows that
\begin{equation}\label{ineq-z-beta-a}
\begin{array}{rcl}
\frac{2|z|}{n^{1/2}} & \le & \frac{2}{n^{1/2}} \| \boldsymbol{\widehat{\beta}} - \boldsymbol{\beta} \|^3_\infty \sum_{i\neq j} | \mu^{\prime\prime}
(\tilde{\beta}_i + \tilde{\beta}_j)|  \\
& \le &  \frac{2}{n^{1/2}} \left\{ \|\bs{\widehat{\beta}}- \bs{\beta} \|_\infty^3 \cdot \sum_{i\neq j } | \mu^{\prime\prime}(\beta_i + \beta_j) |
+ n^2 \|\bs{\widehat{\beta}}- \bs{\beta} \|_\infty^4 \right\} \\
& \lesssim & \frac{ b_n^3 (\log n)^{3/2} }{ n^2 } \sum_{i\neq j } | \mu^{\prime\prime}(\beta_i + \beta_j) |
+ \frac{ b_n^4 (\log n)^2 }{ n^{1/2} } \\
& \lesssim & \frac{ b_n^3 (\log n)^{3/2} }{ n^2 } \sum_{i\neq j } | e^{\beta_i + \beta_j} - 1|
+ \frac{ b_n^4 (\log n)^2 }{ n^{1/2} }
\end{array}
\end{equation}
where the last but one inequality is due to \eqref{ineq-z-beta-a} and the last one inequality is due to \eqref{eq-derivative-mu-various}.
If \eqref{condtition-B} holds, then $|z|/n^{1/2}=o(1)$.
It completes the proof.
\end{proof}

\subsection{Proofs for Theorem 1 (b)}

Let $\widetilde{V}$ denote the Fisher information matrix of $\widetilde{\bs{\beta}}=(\beta_1, \beta_{r+1}, \ldots, \beta_n)^\top$
under the null $H_0: \beta_1 = \ldots= \beta_r$, where
\begin{equation*}
\widetilde{V}=\begin{pmatrix} \tilde{v}_{11} & \bs{\tilde{v}}_{12}^\top \\ \bs{\tilde{v}}_{12} & V_{22} \end{pmatrix},
\end{equation*}
where $V_{22}$ is the lower right $(n-r)\times (n-r)$ block of $V$, $\bs{\tilde{v}}_{12} =
(\tilde{v}_{1,r+1}, \ldots, \bar{v}_{1, n})^\top$, and
\[
\tilde{v}_{11}= 2r(r-1)\frac{ e^{2\beta_1} }{ ( 1 + e^{2\beta_1})^2 }, ~~
\tilde{v}_{1j} = r \frac{ e^{\beta_1 + \beta_j } }{ ( 1 + e^{\beta_1 + \beta_j})^2 },~j=r+1, \ldots, n.
\]
Let $\widetilde{S}=\mathrm{diag}(1/\tilde{v}_{11}, 1/v_{r+1, r+1}, \ldots, 1/v_{nn})$.
With the similar arguments in the proof of Proposition 1 in \cite{Yan:Xu:2013} and in the proof of Theorem 1 in \cite{hillar2012inverses}, we have
\begin{equation}\label{approximate-inverse2-beta}
\|\widetilde{W}:= \widetilde{V}^{-1}-\widetilde{S} \|_{\max} \lesssim \frac{b_n^3}{n^2},
\end{equation}
and
\begin{equation}\label{ineq-inverse2-beta}
\|\widetilde{V}^{-1} \|_\infty \lesssim \frac{b_n}{n},
\end{equation}
respectively.
Recall that $\bs{\widehat{\beta}}^0$ denotes
the MLE of $\bs{\beta} =(\beta_1,  \ldots, \beta_n)^\top$ under the null $H_0: \beta_1=\cdots = \beta_r$,
where $\widehat{\beta}_1^0 = \cdots = \widehat{\beta}_r^0$.
Similar to the proof of Lemma \ref{lemma-consi-beta}, we have:

\begin{lemma}\label{lemma-con-beta-b}
If $b_n=o( (n/\log n)^{1/4} )$, then with probability at least $1-O(n^{-1})$,
\begin{equation*}
\| \bs{\widehat{\beta}}^0 - \bs{\beta}^0 \|_\infty \lesssim b_n \sqrt{ \frac{\log n}{n} } .
\end{equation*}
\end{lemma}

Now, we prove  Theorem \ref{theorem-LRT-beta} (b).

\begin{proof}[Proof of Theorem \ref{theorem-LRT-beta} (b)]
Let  $\mathbf{\tilde{d}}=(\sum_{i=1}^r d_i, d_{r+1},\ldots, d_n)$.
Note that under $H_0$, $\widehat{\beta}_1^0 = \cdots =\widehat{\beta}_r^0$ and $\beta_1=\cdots=\beta_r$.
With the similar arguments in the proof of \eqref{minus-ratio-b}, we have
\begin{equation}\label{likelihood-beta-composite}
\ell_{\beta}(\boldsymbol{\widehat{\beta}}^0)-\ell_{\beta}(\boldsymbol{\beta})
=\frac{1}{2}( \mathbf{\widetilde{d}} - \E \mathbf{\widetilde{d}} )^\top \widetilde{V}^{-1} ( \mathbf{\widetilde{d}} - \E \mathbf{\widetilde{d}} )
-\frac{1}{2}\mathbf{\widetilde{h}}^\top \widetilde{V}^{-1} \mathbf{\widetilde{h}} + \widetilde{z},
\end{equation}
where $\mathbf{\widetilde{h}} = (\tilde{h}_1, \tilde{h}_{r+1}, \ldots, \tilde{h}_n)^\top$, $\tilde{h}_1 = \sum_{i=1}^r h_i$,
\begin{eqnarray*}
h_i & = &  \sum_{j=1,j\neq i}^n \mu^{\prime\prime}( \tilde{\beta}_i + \tilde{\beta}_j) [
\widehat{\beta}_i^0 + \widehat{\beta}_j^0 - (\beta_i + \beta_j)]^2, ~i=1, \ldots, n, \\
\bar{z} & = & \frac{1}{6}\{ \sum_{i=1}^n  (\widehat{\beta}_i^0 -\beta_i)^3 \sum_{j\neq i} \mu^{\prime\prime}( \tilde{\beta}_i + \tilde{\beta}_j)
+ 2 \sum_{i,j=1, j\neq i}^n  (\widehat{\beta}_i^0-\beta_i)^2(\widehat{\beta}_j^0-\beta_j)\mu^{\prime\prime}( \beta_i + \beta_j) \}.
\end{eqnarray*}
In the above equations, $\tilde{\beta}_i$ lies between $\beta_i$ and $\widehat{\beta}_i^0$ for all $i=1, \ldots, n$ and
$\tilde{\beta}_1 = \cdots = \tilde{\beta}_r$.

Note that $r/n \ge \tau>0$ and $\tau$ is a constant.
In view of Lemma \ref{lemma-con-beta-b}, with the similar arguments as in the proof of \eqref{th22}, we have
\begin{equation*}
\frac{ |\mathbf{\bar{h}}^\top \overline{V}^{-1} \mathbf{\bar{h}} | }{\sqrt{n-r}} = o_p(1),~~ \frac{|\bar{z}|}{\sqrt{n-r}}=o_p(1).
\end{equation*}
Now, we evaluate the difference between $( \mathbf{\widetilde{d}} - \E \mathbf{\widetilde{d}} )^\top \widetilde{V}^{-1} ( \mathbf{\widetilde{d}} - \E \mathbf{\widetilde{d}} )$ and $(\dd - \E \dd)^\top V^{-1} ( \dd - \E \dd)$.
By using $\widetilde{S}$ and $S$ to approximate $\widetilde{V}^{-1}$ and $V^{-1}$ respectively,
we have
\begin{eqnarray*}
&&(\dd - \E \dd)^\top V^{-1} ( \dd - \E \dd) -
( \mathbf{\widetilde{d}} - \E \mathbf{\widetilde{d}} )^\top \widetilde{V}^{-1} ( \mathbf{\widetilde{d}} - \E \mathbf{\widetilde{d}} ) \\
& = & \sum_{i=1}^r \frac{ ( d_i- \E d_i )^2 }{ v_{ii} } - \frac{ ( \tilde{d}_1- \E \tilde{d}_1) ^2 }{ \tilde{v}_{11} }
+ ( \mathbf{\widetilde{d}} - \E \mathbf{\widetilde{d}} )^\top \widetilde{W} ( \mathbf{\widetilde{d}} - \E \mathbf{\widetilde{d}} )
+ ( \dd - \E \dd )^\top W ( \dd - \E \dd ).
\end{eqnarray*}
By Lemma \ref{lemma-clt-beta-W}, $( \dd - \E \dd )^\top W ( \dd - \E \dd )=o_p(n)$.
Similar to the proof of Lemma \ref{lemma-clt-beta-W}, we have
\begin{equation*}
( \mathbf{\widetilde{d}} - \E \mathbf{\widetilde{d}} )^\top \widetilde{W} ( \mathbf{\widetilde{d}} - \E \mathbf{\widetilde{d}} )
= o_p(r).
\end{equation*}
Since $\sum_{i=1}^r d_i = 2 \sum_{1\le i<j \le r} a_{ij} + \sum_{i=1}^r \sum_{j=r+1}^n a_{ij}$,
by the central limit theorem for the bounded case (\citet{Loeve:1977}, page 289),
$\bar{v}_{11}^{-1/2}\sum_{i=1}^r (d_i - \E d_i)$ converges in distribution to the standard normal distribution if $\tilde{v}_{11}\to\infty$.
Therefore,
\[
\frac{ [\sum_{i=1}^r \{ d_i-\E(d_i) \}]^2/\tilde{v}_{11} }{ r } = o_p(1).
\]
These arguments show
\[
\frac{1}{2r} \left\{
(\dd - \E \dd)^\top V^{-1} ( \dd - \E \dd) -
( \mathbf{\widetilde{d}} - \E \mathbf{\widetilde{d}} )^\top \widetilde{V}^{-1} ( \mathbf{\widetilde{d}} - \E \mathbf{\widetilde{d}} )
\right\}
= \frac{1}{2r} \sum_{i=1}^r \frac{ ( d_i- \E d_i )^2 }{ v_{ii} }  + o_p(1).
\]
Combining \eqref{minus-ratio-b},
\eqref{th22} and  \eqref{likelihood-beta-composite},  it yields
\begin{equation*}
\frac{ 2 \{ \ell (\boldsymbol{\widehat{\beta}} ) - \ell ( \boldsymbol{\widehat{\beta}}^0 ) \} - r }{ \sqrt{2r}}
= \frac{ \sum_{i=1}^r ( d_i- \E d_i )^2/v_{ii}- r }{\sqrt{2r}}+o_p(1).
\end{equation*}
Similar to Corollary \ref{corollary-weighted-d},
$\sum_{i=1}^r \{ d_i- \E(d_i) \}^2/v_{ii}$ is asymptotically normal distribution with mean $r$ and variance $2r$.
This completes the proof.
\end{proof}

\subsection{Proofs for Theorem \ref{theorem-likelihood-ratio-beta-fixed}}

With some ambiguity of notations,
we let $V_{22}$ denote the Fisher information matrix of $\bs{\beta}^0=(\beta_{r+1}, \ldots, \beta_n)^\top$
under the null $H_0: (\beta_1, \ldots, \beta_r)=(\beta_1^0, \ldots, \beta_r^0)$, where
$V_{22}$ is the lower right $(n-r)\times (n-r)$ block of $V$.
Here, $\beta_1, \ldots, \beta_r$ in the expressions of the elements $v_{ij}$ of $V$ is replaced with $\beta_1^0, \ldots, \beta_r^0$.
Let $S_{22}=\mathrm{diag}( 1/v_{r+1, r+1}, \ldots, 1/v_{nn})$.
With the similar arguments in the proof of Proposition 1 in \cite{Yan:Xu:2013} and in the proof of Theorem 1 in \cite{hillar2012inverses}, we have
\begin{equation}\label{approximate-inverse2-beta2}
\|\widetilde{W}_{22}:= V_{22}^{-1}- S_{22} \|_{\max}  \lesssim \frac{b_n^3}{n^2},
\end{equation}
and
\begin{equation}\label{ineq-inverse2-beta2}
\|V_{22}^{-1} \|_\infty \lesssim \frac{b_n}{n},
\end{equation}
respectively.

An important step in the proof of Theorem \ref{theorem-likelihood-ratio-beta-fixed} is to evaluate
the difference between $\bar{\dd}^\top V^{-1} \bar{\dd}$ and
$\bar{\dd}^\top_2 V^{-1}_{22} \bar{\dd}_2$.
Let $\bar{\dd}_1 = ( \bar{d}_1, \ldots, \bar{d}_r)^\top$ and
\[
V =
\begin{pmatrix} V_{11}  & V_{12} \\
V_{21} & V_{22}
\end{pmatrix}, ~~
W = \begin{pmatrix} W_{11} & W_{12} \\
W_{21} & W_{22},
\end{pmatrix}
\]
where $V_{11}$ and $W_{11}$ have the same $r\times r$ dimension.
By using respectively $S$ to approximate $V^{-1}$ and $S_{22}$ to approximate $V_{22}^{-1}$, we have
\begin{eqnarray}
\label{eq-dVd}
\bar{\dd}^\top V^{-1} \bar{\dd} & = &  \sum_{i=1}^n \frac{ \bar{d}_i^2 }{ v_{ii} }
+ \bar{\dd}_1^\top W_{11} \bar{\dd}_1 + 2 \bar{\dd}_1^\top W_{12} \bar{\dd}_2 + \bar{\dd}_2^\top W_{22} \bar{\dd}_2, \\
\label{eq-dV22d}
\bar{\dd}^\top_2 V^{-1}_{22} \bar{\dd}_2  & = & \sum_{i=r+1}^n \frac{ \bar{d}_i^2 }{ v_{ii} } +  \bar{\dd}_2^\top \widetilde{W}_{22} \bar{\dd}_2.
\end{eqnarray}
To bound $\bar{\dd}_2^\top (W_{22}-\widetilde{W}_{22}) \bar{\dd}_2$, we need to evaluate $\| W_{22}-\widetilde{W}_{22}\|_{\max}$.
Note that
\[
\left[
\begin{pmatrix} S_{11} & \mathbf{0} \\
\mathbf{0} & S_{22}
\end{pmatrix}
+
\begin{pmatrix} W_{11} & W_{12} \\
W_{21} & W_{22}
\end{pmatrix}
\right]
\begin{pmatrix}
V_{11} & V_{12} \\
V_{21} & V_{22}
\end{pmatrix}
=I_{n\times n}
\]
and $S_{22}V_{22} + \widetilde{W}_{22}V_{22} = I_{(n-r)\times (n-r)}$, we have
\[
W_{21}V_{12} + W_{22}V_{22} = \widetilde{W}_{22}V_{22} \Longrightarrow W_{22} - \widetilde{W}_{22}= - V_{22}^{-1}W_{21}V_{12}.
\]
Note that $ V_{22}^{-1} = S_{22} + W_{22}$ and $\| W_{22} \|_{\max} \lesssim b_n^3/n^2 $.
A direct calculation gives that
\begin{eqnarray*}
|(S_{22} W_{21} V_{12})_{ij}| &  = & |\sum_{k=1}^{n-r}\sum_{h=1}^r (S_{22})_{ik} (W_{21})_{k h} (V_{12})_{hj}| \\
& = & | \sum_{h=1}^r \frac{ 1}{v_{i+r,i+r}} (W_{21})_{i h} (V_{12})_{hj} |\\
& \lesssim &  r \cdot \frac{ b_n }{n-1} \cdot \frac{ b_n^3 }{ n^2 } \cdot \frac{1}{4} = O( \frac{ b_n^4 }{ n^3 }),
\end{eqnarray*}
and
\begin{eqnarray*}
|(W_{22} W_{21} V_{12})_{ij}| & = & |\sum_{k=1}^{n-r}\sum_{h=1}^r (W_{22})_{ik} (W_{21})_{k h} (V_{12})_{hj}| \\
& \le & (n-r)r \| W \|_{\max}^2 \frac{1}{4} = O( \frac{ b_n^6 }{ n^3 } ).
\end{eqnarray*}
This shows that
\begin{equation}\label{ineq-W-diff-upper}
\| W_{22} - \widetilde{W}_{22} \|_{\max} \lesssim \frac{ b_n^6 }{ n^3 },
\end{equation}
which is much smaller than $\| W_{22} \|_{\max}$ and $\| \widetilde{W}_{22} \|_{\max}$ themselves.
The order of $n^{-3}$ makes that  $\bar{\dd}_2^\top (W_{22}-\widetilde{W}_{22}) \bar{\dd}_2$ is an asymptotically neglected remainder term.
We have the following lemma.

\begin{lemma}\label{lemma-W-widetilde-d}
If $b_n=o(n^{1/6})$, then
\[
\bar{\dd}_1^\top W_{11} \bar{\dd}_1 = o_p(1),~~ \bar{\dd}_1^\top W_{12} \bar{\dd}_2 = o_p(1),
~~
\bar{\dd}_2^\top ( W_{22} - \widetilde{W}_{22}) \bar{\dd}_2 = o_p(1).
\]
\end{lemma}

Similar to the proof of Lemma \ref{lemma-consi-beta}, we have:

\begin{lemma}\label{lemma-con-beta-c}
Let $\bs{\widehat{\beta}}^0$ be the MLE of $\bs{\beta}^0 =(\beta_{r+1}, \ldots, \beta_n)^\top$ under the null
$H_0: (\beta_1, \ldots, \beta_r)=(\beta_1^0, \ldots, \beta_r^0)$.
If $b_n=o( (n/\log n)^{1/4} )$, then with probability at least $1-O(n^{-1})$,
\begin{equation*}
\| \bs{\widehat{\beta}}^0 - \bs{\beta}^0 \|_\infty  \lesssim b_n \sqrt{ \frac{\log n}{n} } .
\end{equation*}
\end{lemma}

Except for the above lemma, another important step for proving Theorem \ref{theorem-likelihood-ratio-beta-fixed} is to
establish the upper bound of $\max_{i=r+1, \ldots, n} | \widehat{\beta}_i - \widehat{\beta}_i^0 |$, which is stated below.
Note that this error bound has a fast error rate in the magnitude of $n^{-1}$, up to a factor $b_n^3\log n$.
This is much smaller than the error bounds for $\| \bs{\widehat{\beta}}^0 - \bs{\beta}^0 \|_\infty$
and $\| \bs{\widehat{\beta}} - \bs{\beta}^0 \|_\infty$.

\begin{lemma}\label{lemma-hat-beta-diff}
If $b_n^3 = o( n/\log n)$, then with probability at least $1-O(n^{-1})$,
\[
\max_{i=r+1, \ldots, n} | \widehat{\beta}_i - \widehat{\beta}_i^0 | \lesssim \frac{b_n^3 \log n}{n}.
\]
\end{lemma}

Now, we are ready to prove Theorem \ref{theorem-likelihood-ratio-beta-fixed}.
We only state the proof of the first part here. The proof of Theorem \ref{theorem-likelihood-ratio-beta-fixed} (b) is similar and omitted.

\begin{proof}[Proof of Theorem \ref{theorem-likelihood-ratio-beta-fixed} (a)]

Applying a fourth-order Taylor expansion to $\ell_\beta(\widehat{\bs{\beta}}^0 ) $ at point $ \bs{\beta}^0$, it yields
\begin{eqnarray*}
\ell_\beta(\widehat{\bs{\beta}}^0 ) - \ell( \bs{\beta}^0 ) & = &
\underbrace{
\frac{\partial \ell( \bs{\beta}^0 ) }{ \partial \bs{\beta}^\top } ( \widehat{\bs{\beta}}^0 - \bs{\beta}^0 ) +
\frac{1}{2} ( \widehat{\bs{\beta}}^0 - \bs{\beta}^0 )^\top \frac{ \partial^2 \ell( \bs{\beta}^0 ) }{ \partial \bs{\beta} \bs{\beta}^\top } ( \widehat{\bs{\beta}}^0 - \bs{\beta}^0 ) }_{B_1^0} \\
& & + \underbrace{ \frac{1}{6} \sum_{i=1}^n \sum_{j=1}^n \sum_{k=1}^n \frac{ \partial^3 \ell(\bs{\beta}^0)}{ \partial \beta_i \partial \beta_j \partial \beta_k }
( \widehat{\beta}_i^0 - \beta_i^0)( \widehat{\beta}_j^0 - \beta_j^0)( \widehat{\beta}_j^0 - \beta_j^0) }_{B_2^0} \\
&& \underbrace{ \frac{1}{4!} \sum_{t=1}^n \sum_{i=1}^n \sum_{j=1}^n \sum_{k=1}^n \frac{ \partial^4 \ell(\bs{\tilde{\beta}}^0)}{ \partial \beta_t \partial \beta_i \partial \beta_j \partial \beta_k } ( \widehat{\beta}_t^0 - \beta_t^0)
( \widehat{\beta}_i^0 - \beta_i^0)( \widehat{\beta}_j^0 - \beta_j^0)( \widehat{\beta}_j^0 - \beta_j^0) }_{B_3^0},
\end{eqnarray*}
where $\bs{\tilde{\beta}}^0 = \alpha \bs{\beta}^0 + (1-\alpha ) \widehat{\bs{\beta}}$ for some $\alpha\in(0,1)$.
With a similar manner, $\ell_\beta(\widehat{\bs{\beta}} )$ has the following expansion:
\begin{equation}
\ell_\beta(\widehat{\bs{\beta}} ) - \ell( \bs{\beta}^0 ) = B_1 + B_2 + B_3,
\end{equation}
where $B_i$ is the version of $B_i^0$ with $\widehat{\bs{\beta}}^0$ replaced by $\widehat{\bs{\beta}}$.
Therefore,
\[
2[ \ell_\beta(\widehat{\bs{\beta}} ) - \ell_\beta(\widehat{\bs{\beta}}^0 ) ]
= 2( B_1 - B_1^0)  + 2(B_2 - B_2^0) + 2(B_3 - B_3^0).
\]
It is sufficient to demonstrate: (1) $B_1 - B_1^0$ converges in distribution to the Chi-square distribution with $r$ degrees of freedom;
(2) $B_2 - B_2^0$ and $B_3-B_3^0$ are asymptotically neglected remainder terms.
These claims are shown in three steps in turns.

Step 1. We show $B_1 - B_1^0 \stackrel{L}{\longrightarrow}  \chi_r^2$ as $n\to\infty$.
 Let  $\mathbf{d}_2=(d_{r+1}, \ldots, d_n)^\top$, $\mathbf{\hat{d}}=(\hat{d}_1, \ldots, \hat{d}_r, d_{r+1}, \ldots, d_n)^\top$ and
\[
\hat{d}_i = \sum_{j\neq i} \frac{ e^{\widehat{\beta}_i^0 + \widehat{\beta}_j^0 } }{ 1 + e^{\widehat{\beta}_i^0 + \widehat{\beta}_j^0 } },
~~i=1, \ldots, n.
\]
Note that under $H_0$, $( \widehat{\beta}_1, \ldots, \widehat{\beta}_r)=(\beta_1^0, \ldots, \beta_r^0)$.
Let $\boldsymbol{\widehat{\beta}}^0_2= (\widehat{\beta}_{r+1}, \ldots, \widehat{\beta}_n)^\top$ and
$ \boldsymbol{\beta}^0_2 = (\beta_{r+1}^0, \ldots, \beta_n^0)$.
Similar to \eqref{eq-expansion-hatbeta-beta}, we have
\[
\boldsymbol{\widehat{\beta}}^0_2 - \boldsymbol{\beta}^0_2 = V_{22}^{-1}(\dd_2 - \E \dd_2 )  - V_{22}^{-1}\mathbf{h}_2,
\]
where $\mathbf{\widetilde{h}}_2 = (\tilde{h}_{r+1}, \ldots, \tilde{h}_n)^\top$, and
\begin{eqnarray*}
\tilde{h}_i & = &  \sum_{j=1,j\neq i}^n \mu^{\prime\prime}( \tilde{\beta}_i^0 + \tilde{\beta}_j^0) [
\widehat{\beta}_i^0 + \widehat{\beta}_j^0 - (\beta_i^0 + \beta_j^0)]^2, ~i=1, \ldots, n.
\end{eqnarray*}
In the above equation, $\tilde{\beta}_i^0$ lies between $\beta_i^0$ and $\widehat{\beta}_i^0$ for all $i=1, \ldots, n$.
Substituting it into $B_1^0$, it yields
\begin{equation}\label{likelihood-beta-composite}
 B_1^0 = \frac{1}{2}  \mathbf{\bar{d}}_2^\top V_{22}^{-1} \mathbf{\bar{d}}_2 -  \frac{1}{2} \mathbf{\widetilde{h}}_2^\top V_{22}^{-1} \mathbf{\widetilde{h}}_2.
\end{equation}
Note that
\[
B_1 = \frac{1}{2}  \mathbf{\bar{d}}^\top V^{-1} \mathbf{\bar{d}} -  \frac{1}{2} \mathbf{h}^\top V^{-1} \mathbf{h}.
\]
By letting $V^{-1}=S+W$ and $V_{22}^{-1}=S_{22}+ \widetilde{W}_{22}$, we have
\begin{equation}
2(B_1 - B_1^0) = \sum_{i=1}^r \frac{ \bar{d}_i^2 }{v_{ii}}  + \mathbf{\bar{d}}_1^\top W_{11} \mathbf{\bar{d}}_1 +
 2\mathbf{\bar{d}}_1^\top W_{11} \mathbf{\bar{d}}_2 + \mathbf{\bar{d}}_2^\top ( W_{22} - \widetilde{W}_{22}) \mathbf{\bar{d}}_2
+ \mathbf{\widetilde{h}}_2^\top V_{22}^{-1} \mathbf{\widetilde{h}}_2 - \mathbf{h}^\top V^{-1} \mathbf{h}.
\end{equation}
By Lemma \ref{lemma-W-widetilde-d},
\[
\mathbf{\bar{d}}_1^\top W_{11} \mathbf{\bar{d}}_1=o_p(1),~~
\mathbf{\bar{d}}_1^\top W_{11} \mathbf{\bar{d}}_2 =o_p(1),~~ \mathbf{\bar{d}}_2^\top ( W_{22} - \widetilde{W}_{22}) \mathbf{\bar{d}}_2 =o_p(1).
\]
Note that $d_i=\sum_{j\neq i}a_{ij}$ are sums of $n-1$
independent Bernoulli random variables. By the central limit theorem for the bounded case (\cite{Loeve:1977}, p. 289),
$\bar{d}_i/v_{ii}^{1/2}$ converges in distribution to the standard normality if $v_{ii}$ diverges.
Given a fixed $r$, $d_i - \sum_{j=1, j\neq i}^r a_{ij}$, $i=1, \ldots, r$ are independent. Therefore, the vector $(\bar{d}_1/v_{11}^{1/2}, \ldots, \bar{d}_r/v_{rr}^{1/2}$
follows a $r$-dimensional standard normal distribution.
This shows $\sum_{i=1}^r \frac{ \bar{d}_i^2 }{v_{ii}}$ follows a Chi-square distribution with $r$ degrees of freedom.
In the remainder of this step is to show  $\mathbf{\widetilde{h}}_2^\top V_{22}^{-1} \mathbf{\widetilde{h}}_2 - \mathbf{h}^\top V^{-1} \mathbf{h}=o_p(1)$.

Recall the definition of $\mathbf{h}$ in \eqref{eq:definition:h}.
By letting $V^{-1}=S+W$ and $V_{22}^{-1}=S_{22}+W_{22}$, we have
\begin{eqnarray*}
\mathbf{h}^\top V^{-1} \mathbf{h} &  = & \sum_{i=1}^n \frac{ h_i^2 }{ v_{ii} } + \hh_1^\top W_{11} \hh_1
+ 2 \hh_1^\top W_{12} \hh_2 + \hh_2^\top W_{22} \hh_2, \\
\widetilde{\hh}_2^\top V_{22}^{-1} \widetilde{\hh}_2 & = & \sum_{i=r+1}^n \frac{ \tilde{h}_i^2 }{ v_{ii} }
+ \widetilde{\hh}_2^\top \widetilde{W}_{22}^{-1} \widetilde{\hh}_2.
\end{eqnarray*}
By \eqref{ineq-V-S-appro-upper-b} and \eqref{ineq-beta-h}, we have
\begin{eqnarray}
| \hh_1^\top W_{11} \hh_1 | & \le & r^2 \| W_{11} \|_{\max} \| \hh_1 \|_\infty \lesssim r^2 \cdot b_n^2 \log n \cdot \frac{ b_n^3 }{ n^2 },
\\
| \hh_1^\top W_{12} \hh_2 | & \lesssim & r(n-r) \cdot b_n^2\log n \cdot \frac{ b_n^3 }{ n^2 }
\lesssim \frac{ b_n^5 \log n}{ n}.
\end{eqnarray}
To evaluate the bound of $\hh_2^\top W_{22} \hh_2 - \widetilde{\hh}_2^\top \widetilde{W}_{22}^{-1} \widetilde{\hh}_2$,
we divide it into three terms:
\begin{eqnarray*}
& & \hh_2^\top W_{22} \hh_2 - \widetilde{\hh}_2^\top \widetilde{W}_{22}^{-1} \widetilde{\hh}_2 \\
& = & \underbrace{\hh_2^\top W_{22} \hh_2 - \hh_2^\top \widetilde{W}_{22} \hh_2}_{C_1} +
\underbrace{\hh_2^\top \widetilde{W}_{22} \hh_2
- \widetilde{ \hh}_2 \widetilde{W}_{22} \hh_2}_{C_2}
+
\underbrace{\widetilde{\hh}_2^\top W_{22} \hh_2 - \widetilde{\hh}_2 \widetilde{W}_{22} \widetilde{\hh}_2}_{C_3}.
\end{eqnarray*}
The first term $C_1$ is bounded as follows.
By \eqref{ineq-W-diff-upper} and \eqref{ineq-beta-h}, we have
\begin{eqnarray}
\nonumber
| \hh_2^\top ( W_{22} - \widetilde{W}_{22} ) \hh_2 | & \le & (n-r)^2 \| W_{22} - \widetilde{W}_{22} \|_{\max} \| \hh \|_\infty \\
\label{ineq-upper-B1}
& \lesssim & (n-r)^2 \cdot b_n^2 \log n \cdot \frac{ b_n^6 }{ n^3 } \lesssim \frac{ b_n^8}{ n}.
\end{eqnarray}

A key point for bounding $C_2$  is to deriving the error $\| \hh_2 - \widetilde{\hh}_2 \|_\infty$.
For $i=r+1, \ldots, n$, observe that
\begin{eqnarray*}
\tilde{h}_i & = & \left( \sum_{j=1}^r   + \sum_{j=r+1,j\neq i}^n
\right)
\mu^{\prime\prime}( \tilde{\pi}_{ij}^0 ) ( \widehat{\pi}_{ij}^0 - \widehat{\pi}_{ij}^0 )^2
  \\
h_i & = & \left( \sum_{j=1}^r   + \sum_{j=r+1,j\neq i}^n
\right) \mu^{\prime\prime}( \tilde{\pi}_{ij} )
( \widehat{\pi}_{ij} - \widehat{\pi}_{ij}^0 )^2,
\end{eqnarray*}
In order to derive the error bound $| \tilde{h}_i - h_i |$, it is sufficient to bound the following difference:
\begin{eqnarray*}
& & |\mu^{\prime\prime}( \tilde{\pi}_{ij}^0 ) ( \widehat{\pi}_{ij}^0 - \pi_{ij}^0 )^2
-  \mu^{\prime\prime}( \tilde{\pi}_{ij} )
( \widehat{\pi}_{ij} - \pi_{ij}^0 )^2  | \\
& \le & | [\mu^{\prime\prime}( \tilde{\pi}_{ij}^0 )
- \mu^{\prime\prime}(\tilde{\pi}_{ij})]  ( \widehat{\pi}_{ij}^0 - \pi_{ij}^0 )^2 | +
|\mu^{\prime\prime}( \tilde{\pi}_{ij}) [ ( \widehat{\pi}_{ij}^0 - \pi_{ij}^0 )^2 - ( \widehat{\pi}_{ij} - \pi_{ij}^0 )^2]|.
\end{eqnarray*}
where $j=r+1, \ldots, n$, $j\neq i$.
By Lemmas \ref{lemma-consi-beta}, \ref{lemma-con-beta-b}, \ref{lemma-con-beta-c} and \eqref{ineq-mu-deriv-bound}, we have
\begin{eqnarray*}
& & | \mu^{\prime\prime} ( \tilde{\pi}_{ij}^0 ) - \mu^{\prime\prime} ( \tilde{\pi}_{ij} ) | \\
& \le & | \mu^{\prime\prime} ( \tilde{\pi}_{ij}^0 ) - \mu^{\prime\prime}( \hat{\pi}_{ij}^0 ) |
+ | \mu^{\prime\prime}( \hat{\pi}_{ij}^0 )  - \mu^{\prime\prime}( \hat{\pi}_{ij} ) |
+ |  \mu^{\prime\prime}( \hat{\pi}_{ij} ) -  \mu^{\prime\prime}( \tilde{\pi}_{ij} ) | \\
& \le & \frac{1}{4} \| \widehat{\bs{\beta}}^0 - \bs{\beta}^0 \|_\infty
+ \max_{j=r+1, \ldots, n} | \widehat{\beta}_j^0 - \widehat{\beta}_j |
+  \| \widehat{\bs{\beta}} - \bs{\beta}^0 \|_\infty  \\
& \lesssim & b_n \sqrt{\frac{\log n}{n} }.
\end{eqnarray*}
Further, by Lemmas \ref{lemma-consi-beta}, \ref{lemma-con-beta-b}, \ref{lemma-con-beta-c} and \eqref{ineq-mu-deriv-bound}, we have
\begin{eqnarray*}
& &|\mu^{\prime\prime}( \tilde{\pi}_{ij}) [ ( \widehat{\pi}_{ij}^0 - \pi_{ij}^0 )^2 - ( \widehat{\pi}_{ij} - \pi_{ij}^0 )^2]| \\
& \le & \frac{1}{4} \cdot | \widehat{\pi}_{ij}^0 - \widehat{\pi}_{ij}| \cdot ( |\widehat{\pi}_{ij}^0 - \pi_{ij}^0| +
| \widehat{\pi}_{ij} - \pi_{ij}^0 | ) \\
& \lesssim &   \frac{ b_n^4 (\log n)^{3/2} }{ n^{3/2} }.
\end{eqnarray*}
Therefore,
\begin{eqnarray*}
\max_{i=r+1, \ldots, n} | \tilde{h}_i - h_i | & \le & \sum_{i=1}^r | \widehat{\pi}_{ij}^0 - \widehat{\pi}_{ij}^0 )^2 | + \sum_{i=1}^r \mu^{\prime\prime}( \tilde{\pi}_{ij} )
( \widehat{\pi}_{ij} - \widehat{\pi}_{ij}^0 )^2 \\
&& + (n-r) \max_{j=r+1, \ldots, n, j\neq i} |\mu^{\prime\prime}( \tilde{\pi}_{ij}^0 ) ( \widehat{\pi}_{ij}^0 - \pi_{ij}^0 )^2
-  \mu^{\prime\prime}( \tilde{\pi}_{ij} )
( \widehat{\pi}_{ij} - \pi_{ij}^0 )^2  | \\
& \lesssim & \frac{ b_n^3 (\log n)^{3/2} }{ n^{1/2}}.
\end{eqnarray*}
The upper bounds of $C_2$ and $C_3$ are derived as follows.
By \eqref{approximate-inverse2-beta2} and \eqref{ineq-hh2-upper}, we have
\begin{eqnarray}
\nonumber
|C_2|& = & |\hh_2^\top \widetilde{W}_{22} \hh_2
- \widetilde{ \hh}_2 \widetilde{W}_{22} \hh_2 | \\
\nonumber
& \le & ( \hh_2- \widetilde{\hh}_2) ^\top \widetilde{W}_{22} \hh_2 \\
\nonumber
& \le & (n-r)^2 \| \widetilde{W}_{22} \|_{\max} \| \hh_2- \widetilde{\hh}_2 \|_\infty \| \hh_2 \|_\infty \\
\nonumber
& \lesssim & b_n^3 \cdot \frac{ b_n^3 (\log n)^{3/2} }{ n^{1/2} } \cdot b_n^2 \log n  \\
\label{ineq-upper-B2}
& = & \frac{ b_n^8 (\log n)^{5/2} }{ n^{1/2} }.
\end{eqnarray}
By \eqref{ineq-W-diff-upper} and \eqref{ineq-hh2-upper}, we have
\begin{eqnarray}
\nonumber
|C_3|& = & |\widetilde{\hh}_2^\top W_{22} \hh_2 - \widetilde{\hh}_2 \widetilde{W}_{22} \widetilde{\hh}_2 | \\
\nonumber
& \le & | \widetilde{\hh}_2^\top W_{22} \hh_2 - \widetilde{\hh}_2 W_{22} \widetilde{\hh}_2 |
+ |\widetilde{\hh}_2 (W_{22}- \widetilde{W}_{22}) \widetilde{\hh}_2 | \\
\nonumber
& \le & (n-r)^2 \cdot \| W_{22} \|_{\max} \cdot \| \widetilde{\hh}_2 \|_\infty \cdot \| \widetilde{\hh}_2 - \hh_2 \|_\infty
+ (n-r)^2 \cdot \| W_{22} - \widetilde{W}_{22}  \|_{\max} \cdot \| \widetilde{\hh}_2 \|_\infty^2 \\
\nonumber
& \lesssim & n^2 \cdot \frac{ b_n^3 }{n^2 } \cdot \frac{ b_n^3 (\log n)^{3/2} }{ n^{1/2}} \cdot b_n^2 \log n
+ n^2 \cdot \frac{ b_n^6 }{ n^3 } \cdot (b_n^2 \log n) \\
\label{ineq-upper-B3}
& \lesssim & \frac{ b_n^8 (\log n)^{5/2} }{ n^{1/2}}.
\end{eqnarray}
By combining  \eqref{ineq-upper-B1}, \eqref{ineq-upper-B2} and \eqref{ineq-upper-B3}, it yields
\begin{equation}
|\mathbf{h}^\top V^{-1} \mathbf{h} -
\widetilde{\hh}_2^\top V_{22}^{-1} \widetilde{\hh}_2|
\lesssim \frac{ b_n^8 (\log n)^{5/2} }{ n^{1/2} } .
\end{equation}

Step 2. We bound $B_2 - B_2^0$.
Similar to \eqref{cacluation-z}, we have
\begin{equation}
B_2^0 = \frac{1}{6}\left\{ \sum_{i=1}^n  (\widehat{\beta}_i^0 -\beta_i^0)^3 \sum_{j\neq i} \mu^{\prime\prime}( \beta_i^0 + \beta_j^0)
+ 2 \sum_{i=1}^n \sum_{j=1, j\neq i}^n  (\widehat{\beta}_i^0-\beta_i^0)^2(\widehat{\beta}_j^0-\beta_j^0)
\mu^{\prime\prime}( \beta_i^0 + \beta_j^0) \right\}.
\end{equation}
Note that $\widehat{\beta}_i^0=\beta_i^0$, $i=1, \ldots, r$, and
\begin{eqnarray*}
B_2 & = &\frac{1}{6}\{ \sum_{i=1}^n  (\widehat{\beta}_i-\beta_i^0)^3 \sum_{j\neq i} \mu^{\prime\prime}( \beta_i^0 + \beta_j^0 )
+ 2 \sum_{i,j=1, j\neq i}^n  (\widehat{\beta}_i-\beta_i^0)^2(\widehat{\beta}_j-\beta_j)\mu^{\prime\prime}( \beta_i^0 + \beta_j^0 ) \}.
\end{eqnarray*}
For $i,j=r+1, \ldots, n$ and $i\neq j$, we have
\begin{eqnarray*}
| (\widehat{\beta}_i-\beta_i^0)^3 - (\widehat{\beta}_i^0 -\beta_i^0)^3 |
& \le & | \widehat{\beta}_i - \widehat{\beta}_i^0| \cdot [ (\widehat{\beta}_i-\beta_i^0)^2 + |\widehat{\beta}_i-\beta_i^0|
\cdot | \widehat{\beta}_i^0 -\beta_i^0 |+
  (\widehat{\beta}_i^0 -\beta_i^0)^2  ] \\
& \lesssim & \frac{ b_n^3 \log n }{ n } \cdot b_n^2 \left( \frac{\log n}{n} \right),
\end{eqnarray*}
and
\begin{eqnarray*}
&&(\widehat{\beta}_i-\beta_i^0)^2(\widehat{\beta}_j-\beta_j^0) -  (\widehat{\beta}_i^0-\beta_i^0)^2(\widehat{\beta}_j^0-\beta_j^0) \\
& \le & (\widehat{\beta}_i-\beta_i^0)^2(\widehat{\beta}_j-\beta_j^0) - (\widehat{\beta}_i-\beta_i^0)^2(\widehat{\beta}_j^0-\beta_j^0) \\
&& + (\widehat{\beta}_i-\beta_i^0)^2(\widehat{\beta}_j^0-\beta_j^0)-(\widehat{\beta}_i^0-\beta_i^0)^2(\widehat{\beta}_j^0-\beta_j^0) \\
& \le & (\widehat{\beta}_i-\beta_i^0)^2( \widehat{\beta}_j - \widehat{\beta}_j^0 )
+ (\widehat{\beta}_j^0-\beta_j^0)(\widehat{\beta}_i-\beta_i^0) [(\widehat{\beta}_i-\beta_i^0)+(\widehat{\beta}_i^0-\beta_i^0)]\\
& \lesssim & \frac{ b_n^3 \log n }{ n } \cdot b_n^2 \left( \frac{\log n}{n} \right).
\end{eqnarray*}
Therefore,
\begin{equation}\label{ineq-B2-B20}
|B_2 - B_2^0|
\lesssim  r^2 b_n^3 \left( \frac{\log n}{n} \right)^{3/2} + rn b_n^3 \left( \frac{\log n}{n} \right)^{3/2}
+ \frac{ b_n^5 (\log n)^2 }{ n^2 } \sum_{i\neq j} \mu^{\prime\prime}( \beta_i^0 + \beta_j^0 )|.
\end{equation}

Step 3. We bound $B_3 - B_3^0$.
Similar to \eqref{cacluation-z}, we have
\begin{equation}
\begin{array}{rcl}
B_3^0 & = & \frac{1}{4!}\left\{ \sum\limits_{i=1}^n  (\widehat{\beta}_i^0 -\beta_i^0)^3 \sum\limits_{j\neq i} \mu^{\prime\prime\prime}( \tilde{\beta}_i^0 + \tilde{\beta}_j^0)
[(\widehat{\beta}_i^0 -\beta_i^0)+(\widehat{\beta}_j^0 -\beta_j^0)] \right.\\
&& \left.+ 2 \sum\limits_{i=1}^n \sum\limits_{j=1, j\neq i}^n  (\widehat{\beta}_i^0-\beta_i^0)^2(\widehat{\beta}_j^0-\beta_j^0)
\mu^{\prime\prime\prime}( \beta_i^0 + \beta_j^0)[(\widehat{\beta}_i^0 -\beta_i^0)+(\widehat{\beta}_j^0 -\beta_j^0)] \right\}.
\end{array}
\end{equation}
In the above equations, $\tilde{\beta}_i$ lies between $\beta_i$ and $\widehat{\beta}_i^0$ for all $i=1, \ldots, n$.
It is easy to verify
\begin{eqnarray*}
&&(\widehat{\beta}_i -\beta_i^0)^4\mu^{\prime\prime\prime}( \tilde{\beta}_i + \tilde{\beta}_j)-
(\widehat{\beta}_i^0 -\beta_i^0)^4\mu^{\prime\prime\prime}( \tilde{\beta}_i^0 + \tilde{\beta}_j^0) \\
& \le & (\widehat{\beta}_i -\beta_i^0)^4\mu^{\prime\prime\prime}( \tilde{\beta}_i + \tilde{\beta}_j)
- (\widehat{\beta}_i -\beta_i^0)^4 \mu^{\prime\prime\prime}( \tilde{\beta}_i^0 + \tilde{\beta}_j^0) \\
&&+ (\widehat{\beta}_i -\beta_i^0)^4 \mu^{\prime\prime\prime}( \tilde{\beta}_i^0 + \tilde{\beta}_j^0)
- (\widehat{\beta}_i^0 -\beta_i^0)^4\mu^{\prime\prime\prime}( \tilde{\beta}_i^0 + \tilde{\beta}_j^0) \\
& \le & | \widehat{\beta}_i - \widehat{\beta}_i^0| ( |\widehat{\beta}_i- \beta_i^0 |^3 + |\widehat{\beta}_j- \beta_j^0 |^3 )
+ (\widehat{\beta}_i^0 -\beta_i^0)^4 \|\widetilde{\bs{\beta}}_i - \widetilde{\bs{\beta}} \|_\infty \\
& \lesssim & \frac{ b_n^6 (\log n)^{5/2} }{ n^{5/2} }.
\end{eqnarray*}
Similarly, we have
\begin{eqnarray*}
&&
\left | (\widehat{\beta}_i -\beta_i^0)^3\mu^{\prime\prime\prime}( \tilde{\beta}_i + \tilde{\beta}_j)[(\widehat{\beta}_i -\beta_i^0)+(\widehat{\beta}_j -\beta_j^0)] \right.
\\
&&
\left.
-
(\widehat{\beta}_i^0 -\beta_i^0)^3\mu^{\prime\prime\prime}( \tilde{\beta}_i^0 + \tilde{\beta}_j^0)
[(\widehat{\beta}_i^0 -\beta_i^0)+(\widehat{\beta}_j^0 -\beta_j^0)] \right | \\
& \lesssim & \frac{ b_n^6 (\log n)^{5/2} }{ n^{5/2} }. \\
\end{eqnarray*}
Therefore,
\[
|B_3 - B_3^0| \lesssim \frac{ b_n^6 (\log n)^{5/2} }{ n^{1/2} }.
\]
\end{proof}